\documentclass[11pt]{amsart} 
\usepackage{amsthm,amsmath,amssymb,latexsym}
\usepackage{hyperref}
\hypersetup{backref, pdfpagemode=FullScreen, colorlinks=true}
\begin{document}

\oddsidemargin=18pt \evensidemargin=18pt
\newtheorem{theorem}{Theorem}[section]
\newtheorem{definition}[theorem]{Definition}
\newtheorem{proposition}[theorem]{Proposition}
\newtheorem{lemma}[theorem]{Lemma}
\newtheorem{remark}[theorem]{Remark}
\newtheorem{corollary}[theorem]{Corollary}
\newtheorem{question}{Question}
\newtheorem{fact}{Fact}[section]
\newtheorem{claim}[theorem]{Claim}

\newtheorem{notation}{Notation}[section]
\newcommand\rem{\begin{remark}\upshape}
\newcommand\erem{\end{remark}}
\newcommand\ft{\begin{fact}\upshape}
\newcommand\eft{\end{fact}}
\newcommand\nota{\begin{notation}\upshape}
\newcommand\enota{\end{notation}}
\newcommand\dfn{\begin{definition}\upshape}
\newcommand\edfn{\end{definition}}
\newcommand\cor{\begin{corollary}}
\newcommand\ecor{\end{corollary}}
\newcommand\cl{\begin{claim}}
\newcommand\ecl{\end{claim}}
\newcommand\thm{\begin{theorem}}
\newcommand\ethm{\end{theorem}}
\newcommand\prop{\begin{proposition}}
\newcommand\eprop{\end{proposition}}
\newcommand\lem{\begin{lemma}}
\newcommand\elem{\end{lemma}}
\providecommand\qed{\hfill$\quad\Box$}
\newcommand\pr{{Proof:\;}}
\newcommand\prcl{\par\noindent{\em Proof of Claim: }}
\newcommand\dom{{\text{dom}}}
\newcommand\Pn{{P_n}}
\newcommand\ord{{ord}}
\newcommand\rk{{rk}}
\newcommand\G{{G^{\n}}}
\newcommand\car{{\rm{char}}}
\newcommand\M{{{\mathcal M}}}
\newcommand\N{{{\mathcal N}}}
\newcommand\K{{{\mathcal K}}}
\newcommand\D{{{\mathcal D}}}
\newcommand\A{{{\mathcal A}}}
\newcommand\Ps{{{\mathcal P}}}
\newcommand\E{{\mathcal E}}
\newcommand\X{{{\bold{X}}}}
\newcommand\Y{{{\bold{Y}}}}
\renewcommand\O{{{\mathcal O}}}
\renewcommand\L{{{\mathcal L}}}
\newcommand\Th{{\text{Th}}}
\newcommand\IZ{{\mathbb Z}}
\newcommand\IQ{{\mathbb Q}}
\newcommand\IR{{\mathbb R}}
\newcommand\IN{{\mathbb N}}
\newcommand\IC{{\mathbb C}}
\newcommand\F{{\mathbb F}}
\newcommand\Fe{{\mathcal F}}
\newcommand\Se{{\mathcal S}}
\newcommand\V{{\mathcal V}}
\newcommand\W{{\mathcal W}}
\newcommand\T{{\mathcal T}}
\newcommand\si{{\sigma}}
\newcommand\n{{\nabla}}
\newcommand\C{{\mathcal C}}
\newcommand\Lr{{{\mathcal L}_{\text{rings}}}}
\newcommand\Lrd{{{\mathcal L}^*_{\text{rings}}}}
\newcommand\B{{\mathcal B}}
\newcommand\La{{\mathcal L}}
\def\U{{ \mathfrak U}}
\def\B{{\mathfrak B}}
\def\I{{\mathcal I}}
\let\le=\leqslant
\let\ge=\geqslant
\let\subset=\subseteq
\let\supset=\supseteq
\title[Definable groups in topological differential fields]{Definable groups in topological differential fields}

\author{Fran\c coise Point}
\address{\hskip-\parindent
  Fran\c coise Point\\
  D\'epartement de Math\'ematique\\
  Universit\'e de Mons, Le Pentagone\\
  20, Place du Parc, B-7000 Mons, Belgium}
\email{point@math.univ-paris-diderot.fr}
\thanks{Research Director at the "Fonds de la Recherche
  Scientifique-FNRS".}

\begin{abstract} 
For certain theories of existentially closed topological differential fields, we show that there is a strong relationship between $\L\cup\{D\}$-definable sets and their $\L$-reducts, where $\L$ is a relational  expansion of the field language and $D$ a symbol for a derivation. This enables us to associate with an $\L\cup\{D\}$-definable group in models of such theories, a local $\L$-definable group. As a byproduct, we show that in closed ordered differential fields, one has the descending chain condition on centralisers. 
\end{abstract}
\maketitle
{\it MSC:} 03C10, 03C60, 12J,12H05, 22A05\\
\par{\it Keywords: NIP theories, derivation, topological fields, definable groups}
\section{Introduction}
\par Let $\K$ be an $\L$-structure expanding a field of characteristic $0$, and assume that $\K$ is endowed with a non-discrete definable field topology.  
We assume that $\L$ is a relational expansion of the ring language possibly with additional constants (with further assumptions on the relation symbols, as defined in \cite{GP}, see section 2) and that $\K$ is a model of a universal  $\La\cup \{^{-1}\}$-theory $T$ which admits a model-completion $T_{c}$. Furthermore we suppose that the class of models of $T_{c}$ satisfies Hypothesis
$(I)$ which generalizes the notion of {\it large fields} \cite{Pop} $(\star)$. 
Examples of such theories $T_c$ are the theories of algebraically closed valued fields, real-closed fields, real-closed valued fields and p-adically closed fields. In particular these theories are dp-minimal and for the first three respectively $C$-minimal, o-minimal, weakly o-minimal.
\par Note that W. Johnson showed that if $\K$ is an expansion of a field, $K$ infinite, and $dp$-minimal (a stronger condition than NIP) but not strongly minimal, then $\K$ can be endowed with a non-discrete definable field topology, namely $\K$ has a uniformly definable basis of neighbourhoods of zero compatible with the field operations \cite[Theorem 9.1.3]{JW}. Furthermore the topology is induced either by a non-trivial valuation or an absolute value. Under the stronger hypothesis of being dp-small, such field $K$ is either algebraically closed or real-closed \cite[Theorem 9.1.7]{JW}.
\par We consider the {\it generic} expansion of $\K$ with a derivation $D$, namely we put no a priori continuity assumptions on $D$. In certain classes of existentially closed  differential expansions $\K_D$, we want to relate the properties of the definable subsets and then definable subgroups of $\K_D$ with those in the reduct $\K$.
\par Let $\L_D:=\L\cup \{^{-1}\}\cup\{D\}$. Under the above assumption $(\star)$ on $T$, one can show \cite{GP} that the expansion of $T$
to the $\L_D$-theory $T_{D}$ consisting of $T$ together with the axioms expressing that $D$ is a derivation, admits a model-completion $T_{c,D}^*$. An axiomatization of $T_{c,D}^*$ can be obtained as follows: to the theory $T_{D}\cup T_{c}$, one adds a scheme of axioms (DL), which expresses that each differential polynomial has a zero close to a (regular) zero of its associated algebraic polynomial. 
\par This scheme (DL) is related to the axiom scheme (UC) introduced by M. Tressl in the framework of large fields, and in the case of existentially closed ordered differential fields, it provides an axiomatization alternative to the one given by M. Singer (denoted by CODF). 
A number of properties of CODF generalize to this larger setting.
\par Using the same strategy as in \cite{D89}, one defines a fibered dimension function on the definable subsets in models of $T_{c,D}^*$ \cite{GP2}, proving that $T_{c,D}^*$ has the equational boundedness property \cite[Corollary 3.10]{GP2}.
For models of $CODF$, a fibered dimension function has been already introduced but using a cell-decomposition theorem \cite{BMR}. One can show that both dimension functions coincide, checking it on definable subsets of the domain.
\par Using that $T_{c}$ admits quantifier elimination, an easy observation is that whenever $T_{c}$ has $NIP$, then $T_{c,D}^*$ has $NIP$ \cite{GP}. This implies in particular that the definable subsets in a model of $T_{c,D}^*$ have a $VC$-dimension and in \cite{P2}, we related the VC-density of definable sets in $T_c$ and in $T_{c,D}^*$. (Note that $T_c$ dp-small does not imply that $T_{c,D}^*$ dp-small, e.g. CODF is not even strongly dependent \cite[Lemma 4.5.9]{Br}.)
\par Here under the additional assumption that $T_c$ has finite Skolem functions and satisfies the local continuity property (see section 3), we will first prove that there is a strong relationship between $\La_D$-definable sets and $\La$-definable sets in models of $T_{c,D}^*$. 
As a consequence we will show that $T_{c,D}^*$ eliminates the quantifier $\exists^{\infty}$, which will imply for instance that an infinite externally definable set contains an infinite internally definable set \cite[Corollary 9.2.11]{JW}. This will also entail that the open core of a model of $T_{c,D}^*$ is dp-minimal, respectively o-minimal, $p$-minimal, weakly-o-minimal whenever the theory $T_c$ has that property (see section 3).
\par Again this was first done in models of CODF and as an application, one can show that CODF admits elimination of imaginaries, using the corresponding result for RCF \cite{P1}. Note also that any model of CODF is definably complete (section 3), which entails a number of properties on definable continuous (respectively monotone) functions\cite{servi}, \cite{FH}.  
\par Using again the strong relationship between $\La_D$-definable sets and $\La$-definable sets in models of $T_{c,D}^*$, together with a strong density property of definable types in CODF, one can give another proof of elimination of imaginaries \cite{BCP} using a criterium due to E. Hrushovski \cite{H14} which relies on a precise description of definable types.
\par Then, we will associate with any $\L_{D}$-definable group $G$ in a model of $T_{c,D}^*$ an $\L$-definable local group $G^*_{local}$, where the differential points coming from $G$ are dense. If $G$ has a definable subgroup $G_{0}$, we show that $G_{0}$ corresponds to a definable sublocal subgroup of $G^*_{local}$.
Recall that local groups arise for instance as Lie groups of transformations \cite[chapter 1]{MZ}; I. Goldbring proved, using methods from  nonstandard analysis, that a local group which is locally euclidean is locally isomorphic to a Lie group \cite{G}. 
\par Then we apply our result to the special case when $T_{c}$ is the theory of real-closed fields and show that $\L_{D}$-definable groups in closed ordered fields have the descending chain condition on centralisers. 
\par Finally we consider type-definable subgroups. Assuming that $T_{c}$ is NIP, given a definable group $G$, we show how to associate with $G^{00}$, the smallest type-definable subgroup of $G$ of bounded index, an equivalence relation $\sim^{00}$ on a large subset of $G^*_{local}$ such that each equivalence class has a representative in $G/G^{00}$.  

\section{Preliminaries}
\par Let $\La$ be a relational expansion of the language of rings, $\La:=\La_{\text{rings}} \cup \{R_i; i\in
I\}\cup\{c_j;j\in J\}$ where the $c_j$'s are constants and the $R_i$ 
are $n_{i}$-ary predicates. Let $\K$ be an $\La\cup \{^{-1}\}$-structure such that its restriction to $\La_{\text{rings}}\cup \{^{-1}\}$ is field. Let $\tau$ be a topology on $K$ and assume that $\langle K,\tau\rangle$ is a topological $\La$-field as defined in \cite{GP}. In particular, $K$ will always be a field of characteristic $0$ and 
every relation $R_{i}$ (respectively its complement $\lnot R_{i}$), with $i\in I$, is interpreted
in $K$, as the union of an open set $O_{R_{i}}$
(respectively $O_{\lnot R_i}$) and an algebraic subset $\{\bar x\in
K^{n_{i}}:\bigwedge_k r_{i,k}(\bar x)=0\}$ of $K^{n_i}$ (respectively
$\{\bar x\in K^{n_{i}}:\bigwedge_{l} s_{i,l}(\bar x)=0\}$ of $K^{n_{i}}$), where $r_{i,k},\;s_{i,l}\in K[X_1,\cdots,X_{n_i}]$.)
Examples of such structures are given in \cite[Section 2]{GP}.
\par Let $\L_D:=\La \cup\{^{-1}, D\}$ where $D$ is a new unary function symbol.
We consider {\it differential topological $\La$-fields}, namely expansions of topological $\La$-fields to $\L_D$-structures
which satisfy the following two axioms expressing that $D$ is a derivation:
\begin{equation*}
\forall a\,\forall b \;\;D(a+b)=D(a)+D(b)\text{, }
\forall a\,\forall b\;\;\; D(a.b)=a.D(b)+D(a).b\text{.}
\end{equation*}
Unlike other classes of differential topological fields with a well-understood model theory, for instance the $D$-valued fields studied by T. Scanlon \cite{Sc}, or the field of transseries studied by M. Aschenbrenner, L. van den Dries and J. van der Hoeven \cite{ADH}, we impose a priori no continuity assumptions on $D$.
\medskip
\nota\label{dfninitial} Let $\bar a\in K^n$, we denote its length by $\ell(\bar x)$. Set $D(\bar a):=(D(a_{1}),D(a_{2}),\cdots,D(a_{n}))$ and denote the tuple $(a_{1},D(a_{1}),\cdots, D^m(a_{1}))$ by $a_{1}^{\n_{m}}$ and by $\bar a^{\n_{m}}$ the tuple $(a_{1}^{\n_{m}},\cdots,a_{n}^{\n_{m}})$.
\par Let $K\{X_{1},\cdots,X_{n}\}$ be the  ring of differential polynomials over $K$ in $n$ differential indeterminates $X_{1},\cdots, X_{n}$ over $K$, namely it is the ordinary polynomial ring $K[X_{i}^{(j)}; 1\leq i\leq n, j\in \IN]$, with by convention $X_{i}^{(0)}:=X_{i}$. One can extend the derivation $D$ of $K$ to this ring by setting $D(X_{i}):=X_{i}^{(1)}$ and using the additivity and the Leibnitz rule. 
\par Set $\X:=X_{1},\cdots, X_{n}$ and for $\bold{k}:=(k_{1},\cdots,k_{n})\in \IN^{n}$, let $\vert \bold{k} \vert:=\sum_{i=1}^n k_{i}$ and $\X^{\bold{k}}:=X_{1}^{k_{1}}\cdots X_{n}^{k_{n}}$. We denote the usual partial derivative by  $\partial_{\bold{k}}:=(\frac{\partial}{\partial X_{1}})^{k_{1}}\cdots (\frac{\partial}{\partial X_{n}})^{k_{n}}$.
We will use the lexicographic ordering on $\IN^{n}$.
\par Let $f(\X)\in K\{\X\}\setminus K$ and suppose that $f$ is of order $m$, then we write 
$f(\X)$ as $f^*(X_{1},\ldots ,X_{1}^{(m)},\cdots,X_{n},\ldots ,X_{n}^{(m)})$ for some ordinary polynomial 
$f^*(X_{1},\cdots,X_{n.(m+1)})\in K[X_{1},\cdots,X_{n.(m+1)}].$ 
We will make the following abuse of notation: if $\bar b\in K^n$, then $f^*(\bar b)$ means that we evaluate the polynomial $f^*$ at the tuple $\bar b^{\n_{m}}$. (When there is no danger of confusion, we will simply denote it by $\bar b^{\n}$).
\par If $n=1$, recall that the separant $s_{f}$ of $f$ is defined as $s_{f}:=\frac{\partial f }{\partial X_{1}^{(m)}}$.
\par If $g(\X)$ is some ordinary polynomial in $K[X_{1},\cdots,X_{n}]$ of total degree $k$ (i.e. $k=max \{\vert \ell\vert:\;{\X}^{\ell}$ occurs (with non-zero coefficients) in $g(\X)\}$, we will use the Taylor expansion of $g$:
$$g(\X+\Y)=g(\X)+\sum_{\vert \bold{\ell}\vert=1}^{k} g_{\ell}(\X).\Y^{\bold{\ell}},$$ 
where $\ell!g_{\ell}:=\partial_{\bold{\ell}} g$. 
\enota
\par Let $R:=K[\X]$. Recall that $R[X_{n+1}]$ satisfies the generalized Euclidean algorithm. Namely, given $f(X_{n+1}), g(X_{n+1})\neq 0$ two polynomials in $R[X_{n+1}]$ and letting $b\in R$ be the leading coefficient of $g(X_{n+1})$ (viewed as a polynomial in $X_{n+1}$).
Then there exists $d\in \IN$ and polynomials $q(X_{n+1}), r(X_{n+1})\in R[X_{n+1}]$ such that $b^d.f(X_{n+1})=q(X_{n+1}).g(X_{n+1})+r(X_{n+1})$ with $deg_{X_{n+1}}(r(X_{n+1}))<deg_{X_{n+1}}(g(X_{n+1}))$ \cite[Theorem 2.14]{J}.
\nota We will denote the projection maps between cartesian products of $K$ as follows: $\pi^n_{(i_{1},\cdots,i_{k})}:K^n\rightarrow K^k:(x_{1},\cdots,x_{n})\rightarrow (x_{i_{1}},\cdots,x_{i_{k}}),$ $1\leq k\leq n$, $i_{1}<\cdots<i_{k}\leq n$. In the case where $(i_{1},\cdots,i_{k})=(1,\cdots,k)$, we denote the corresponding map by $\pi_{k}$ (not to be confused with the projection on the $k^{th}$ coordinate denoted by $\pi_{(k)}$).
For ease of notation, we will also denote $\pi^n_{(i_{1},\cdots,i_{k})}$ by $\pi^n_{\bar c}$, where $\bar c$ is a tuple of $0$'s and $1$'s with the $1$'s at the $i_{1},\cdots,i_{k}$ positions and $0$'s elsewhere. Also given two tuples $\bar c_1,\;\bar c_2$ of $0$'s and $1$'s, we denote the corresponding projection $\pi_{\bar c_1^{\frown}\bar c_2}$ by $\pi_{\bar c_1\bar c_2}$.
\enota
\par We will consider a class $\mathcal{C}$ of topological
$\La$-fields $K$ of {\it characteristic zero} 
which have a uniformly definable basis of neighbourhoods of zero,
 namely  there is a formula $\varrho(x,\bar y)$ such that the set of subsets of the form $\varrho(K,\bar a):=\{x\in K:\;K\models \varrho(x,\bar a)\}$ with $\bar a\subset K$ 
can be chosen as a basis $\V$ of neighbourhoods of $0$ in $K$ \cite{P87}. 
\par From now on we will assume that $\C$ is the class of models of a universal $\La\cup \{^{-1}\}$-theory $T$ which admits a model-completion $T_{c}$.
We further assume that the class of models of $T_c$ satisfies Hypothesis
$(I)$ which generalizes in our topological setting of the notion of {\it large fields} \cite{Pop}. 
\par Let $T_D$ (respectively $T_{c,D}$) be the $\La_{D}$-theory $T$
(respectively $T_c$) together with the axioms 
stating that $D$ is a derivation.
In \cite{GP}, we showed that under the assumption that the models of $T_{c}$ satisfies Hypothesis $(I)$, that any model of $T_c$ embeds in an element of $T$ satisfying the scheme $(DL)$, namely that  if for 
each differential polynomial
 $f(X)=f^*(X,X^{(1)},\ldots ,X^{(n)})\in K\{X\}$,
for every $W\in \V$, 

$\displaystyle (\exists \alpha_0,\ldots ,\alpha_n\in K)
(f^*(\alpha_0,\ldots ,\alpha_n)=0 \wedge s_f^*(\alpha_0,\ldots ,\alpha_n)\ne 0 )
\Rightarrow$\\
$\displaystyle \Big((\exists z)\big(f(z)=0\wedge s_f(z)\ne 0\wedge$
$\displaystyle \bigwedge_{i=0}^n (
z^{(i)}-\alpha_i\in W)\big)\Big)$.

\par In \cite{GP}, the scheme (DL) was stated only for polynomials with coefficients in a certain neighbourhood of zero and then it was shown that it is equivalent to the present form \cite[Proposition 3.14]{GP}.
\par Let $T_{c,D}^*$ be the $\La_{D}$-theory consisting of $T_{c,D}$ together with the scheme $(DL)$.

\ft\label{qe}\rm{(\cite[Theorem 4.1, Corollary 4.3]{GP})}
Under the above hypotheses on $T$ and $T_c$, 
we have that the theory $T_{c,D}^*$ is the model-completion of $T_D$ and that the theory $T_{c,D}^*$ admits quantifier elimination. Moreover, if $T_{c}$ and $T_{c,D}^*$ are complete and if $T_{c}$ has $NIP$, then $T_{c,D}^*$ has $NIP$. \qed
\eft
\nota\label{star} Let $\phi(x_{1},\cdots,x_{n})$ be an open $\La_{D}$-formula, for each $x_{i}$, $1\leq i\leq n$, let $m_{i}$ be the maximal natural number $m$ such that $D^m(x_{i})$ occurs in an atomic subformula. We will also denote $D^m(x_{i})$ by $x_{i}^{(m)}$.
Then, we denote by $x^*$ the tuple $(x_{i,j})_{i=1,j=0}^{n,m_{i}}$ and by $\phi^*(x^*)$ the $\La$-formula we obtain from $\phi$ by replacing each $x_{i}^{(j)}$ by $x_{i,j}$. We will allow ourselves to denote that formula by $\phi^*$.
\par Let $N:=\sum_{i=1}^n m_{i}+n$ and if $S$ is the subset of $K^n$ defined by $\phi$, we denote by $S^{alg}$ the subset of $K^N$ definable by $\phi^*$. 
\par Note that given a definable set $S:=\phi(K)$, the {\it operation} of taking $S^{alg}$ depends on $\phi$. However, if $A\subset B$ are two definable subsets of $K^n$, then by possibly taking the direct product with a power of $K$, one may assume that $A^{alg}\subset B^{alg}$. 
\par If we wish to stress that the parameters are taken from $K$, we speak of $\La_{D,K}$ (or $\L_{K}$) definable subsets. Denote by $Def(\K)$ the set of  all definable subsets (in any cartesian powers of $K$) of $\K$, with parameters from $K$. 
\enota

\dfn A definable subset $S\in Def(\K)$ has a non-empty $*$-interior if the corresponding $\L_{K}$-definable set $S^{alg}$ has a non-empty interior.
We will denote the interior of a set $A$ by $Int(A)$ and the relative interior of $A\subset B$ in $B$ by $Int_{B}(A)$. 
\edfn

\par Using a notion of independence built in from the algebra of all terms (t-independence), following \cite{D89}, \cite{Marc}, one can define on $Def(\K)$, when $\K\models T_{c,D}^*$ a fibered dimension function (t-dim$_K$)  \cite{GP2}. In order to show that this dimension is fibered, one uses the closure operator $cl$ which is defined by:  
$a\in cl(A)$ if and only if there is a differential polynomial $q\in K \langle A\rangle\{X\}\setminus\{0\}$ such that $q(a)=0$. 
Then associated with this closure operator, one has a notion of independence and dimension for tuples of elements. Let $\K^*$ be a $\vert K\vert^+$-saturated extension of $\K$.
Define for $\bar a=(a_1,\ldots,a_n)$ a tuple in $\K^*$,
cl-$\dim(\bar a)$ := $\max \{ |B| : B\subseteq K\langle \bar a\rangle, B \;{\rm is\; cl-independent} \}.$  

Then one constructs a function $t$-dim on $K$-definable sets $\varphi(K)$ and shows that \cite[Lemma 2.11]{GP2}
$t$-dim$_K(\varphi(K))=\max\{cl$-dim$(\bar c):\;c\in \varphi(\K^*)\}$ $\;\;(\star).$
One then proves that any model of $T_{c,D}^*$ is equationally bounded, which entails that $t$-dim$_K$ is a fibered dimension function \cite[Corollary 3.10]{GP2} (and see below). This dimension function has been further investigated in \cite{BCP} and one can check that for $\bar a$ a tuple of elements in $K^*$, that $cl$-dim($\bar a$)=inf$\{t$-dim$(\varphi(K^*)):\; \varphi(\bar x)\in tp(\bar a/K)\}$. Note that there are infinite subsets of $t$-dimension $0$.
\medskip

\par A {\it generic} element in $\phi(K)$ is a tuple with the property that cl-$dim(\bar a)$=$t$-dim$(\varphi(K^*))$. We will denote by $\L_{D}$-generic, in order to distinguish it from $\L$-generic which corresponds to the model-theoretic closure operator {\it $acl^{\L}$} (or $acl^{\L_{K}}$) occurring in models of $T_{c}$. 

\ft\label{cor delta dim}\cite[section 3]{GP2}
Let $\mathcal K\models T_{c,D}^*$, then  
there exists a fibered dimension function $t$-dim on $Def(\K)$ with the following properties:
\begin{enumerate}
\item Let $S\in Def(\K)$ 
either with non empty interior or with non empty $*$-interior, then $S$ is not included in a zero-set of a non trivial (differential) polynomial and so, in particular, $t$-dim$(S)\neq 0$.
\item Let $A\subset B\in Def(\K)$ and assume that $A^{alg}\subset B^{alg}$. Let $k\geq 1$, then $t$-dim$(A)<t$-dim$(B)=k$ iff there exists a projection $\pi_{k}$ in $K^k$ such that the relative $*$-interior of $\pi_{k}(A)$ in $\pi_{k}(B)$ is empty, whereas the $*$-interior of $\pi_{k}(B)$ is non-empty. 
\item Let $X\in Def(\K)$ and let $f$ be a definable bijection. Then $t$-dim$(X)\geq k$ if and only if $t$-dim$(f(X))\geq k$, $k\in \IN$.
\end{enumerate}
\eft
\rem In the special case of the theory CODF of closed ordered differential fields, Brihaye, Michaux and Rivi\`ere defined a dimension function ($\delta$-dim) on definable subsets in models $\M$ of $CODF$, proving first a cell decomposition theorem \cite[Theorem 4.9]{BMR}. They proved that the dimension function $\delta$-dim is fibered \cite[Theorem 5.19]{BMR}. 
Using the Fact above, one obtains
a possibly other fibered dimension function on definable subsets in models of $CODF$ \cite[Corollary 3.10]{GP}. It is easily noted that these two dimension functions coincide; either one 
notes that on definable subsets of $M$, they do coincide, or one uses the characterisation $(\star)$ \cite[Corollary 5.27]{BMR}, \cite[Lemma 2.11]{GP}.
\erem
\dfn \label{dlarge} Let $\M$ be an $\L$-structure endowed with a fibered dimension function $\dim$. Then a definable subset $S$ of $M^n$ is {\it large} if $\dim(M^n\setminus S)<\dim(S)$ \cite[Definition 1.11]{P88}. 
\edfn
\par We will use repeatedly the following result which can be found in \cite{P88}, in the setting of o-minimal theories but it only uses the notion of a {\it fibered dimension}.
\ft \label{def-large} {\rm \cite[Proposition 1.13, Remark 1.14]{P88}} Let $\M$ be an $\L$-structure endowed with a fibered dimension function. Let $A$ be an $\L_{B}$-definable subset of $\M$, with $B\subset M$ and let $\phi(\bar x;\bar y)$ be an $\L$-formula. Then $\{\bar d:\;\phi(M;\bar d)\cap A$ is large in $A\}$ is $B$-definable.
Moreover it is equal to $\{\bar d:$ for every generic point $\bar u$ of $A$ over $\bar d,$ $\M\models \phi(\bar u;\bar d)\}$.
\eft
\par Using Fact \ref{cor delta dim}, we may deduce the following properties of a definable group $G$ in a model $\K$ of $T_{c,D}^*$.
\ft \label{generic} {\rm \cite{P88}} Any element of $G$ is the product of two $\L_{D}$-generics and given a $\L_{D}$-generic $a$ of $G$ over $b$, then $b.a$ is a $\L_{D}$-generic in $G$.
\eft
\dfn {\rm \cite{P88}} Let $G$ be a definable group of $\K\models T_{c,D}^*$. A {\it generic} definable subset of $G$ is a subset of $G$ such that finitely many translates cover $G$. 
\par Assume that $G$ is defined over $K_{0}$, where $\K_{0}\prec \K$ and let $X$ be a definable or type-definable subset of $G$. Then $X$ is {\it finitely satisfiable} if $X\cap M_{0}^n\neq \emptyset$.
\edfn
Note that it entails that a generic set has the same $t$-dimension than $G$ and so it contains a $\L_{D}$-generic point of $G$.
\ft {\rm \cite[Lemma 2.4]{P88}} Let $G$ be a definable group of $\K\models T_{c,D}^*$. Let $X$ be a definable large subset of $G$, then $X$ is finitely satisfiable and generic. 
\eft
\medskip

\nota \label{sharp} Given a tuple $\bar b\in K^{n}$, we will say that it is a {\it differential tuple}, if there exists $b\in K$ such that $\bar b=b^{\nabla_{n-1}}$, $n\in \IN^{>0}.$
\par Given two tuples $\bar b,\; \bar c$ with $\bar b\in K^{n+1}$ and $\bar c\in K^{n+m+1}$, $1\leq n$, $1\leq m$, we will say that $\bar c$ is a {\it differential prolongation} of $\bar b$ if there exists $b\in K$ such that $\bar b=b^{\nabla_{n}}$ and $\bar c=b^{\nabla_{n+m}}$.
\par Let $A\subset K^n$ be a set of differential tuples and $B$ a subset of $K^{n+m}$ for some $n, m\in \IN^*$. Then $B$ is a {\it differential prolongation} of $A$ if every element of $B$ is a differential prolongation of some element of $A$ and every element of $A$ has a differential prolongation in $B$.
\par Given an $\L_{D}$-definable subset $A=\phi(K)\subset K^m$ with $\phi(x_{1},\cdots,x_{m})$ a $\L_{D}$-formula, and the corresponding $\L$-definable set $A^{alg}=\phi^*(K)\subset M^{m.(n+1)}$, $n\in \IN$, with $\phi^*$ a $\L$-formula, we want to single out inside this $\L$-definable set, the elements of the form $\bar a^{\n}$, with $\bar a\in K^m$.  We will denote the subset of $A^{alg}$ of such elements by $A^{\n}$. 
\enota
\par Finally, let us adapt the definition of an {\it open core} \cite[Definition, page 1372]{DMS} to our setting.
\dfn\label{core}  Let $\K$ be a differential topological $\La$-field. The open core of $\K$ is the relational structure denoted by $\K^{o}$ consisting of the domain of $\K$ together with a relation symbol for each non-empty element of $Def(\K)$ which is open.
\edfn
Note that the reduct of $\K$ to an $\La$-structure is a model of $T_{c}$, then it is interdefinable with $\K^{o}$ \cite[section 4, 4.1(1)]{DMS}. To prove this, one uses that $T_{c}$ admits quantifier elimination and a former result of R. Dougherty and C. Miller who showed that in a first-order topological structure, a definable set which is a boolean combination of open sets is a boolean combination of definable open sets \cite{DM}.

\section{Definable sets in the Mathews' setting}
\par As before, let $\L$ be an expansion of the language of rings with extra constants and relation symbols (but no new function symbols).
\par L. Mathews introduced a class of first-order topological $\L$-structures for which he proved a cell decomposition theorem. Then he applied his result to define a fibered dimension function on the definable subsets of elements of that class \cite{M}. Instances of such classes are the class of o-minimal structures, p-adically closed fields, real-closed rings and real-closed valued fields. 
\par We will put further hypothesis on the theory $T_c$ to be in the same setting, but let us first recall the definition of the cell decomposition property (CDP) for a topological structure.
\dfn \cite[section 6]{M} Let $\A$ be a first-order topological $\L$-structure. A {\it cell} consists either of a point in $A^n$ and is of dimension $0$, or is a definable subset $X$ of $A^n$ such that there exist $1\leq k\leq n$ and $k$ positive integers $1\leq i_{1}<\cdots<i_{k}\leq n$ such that the projection map $\pi^n_{(i_{1},\cdots,i_{k})}$ is a definable homeomorphism and $\pi^n_{(i_{1},\cdots,i_{k})}(X)$ is an open subset of $A^k$. The {\it dimension} of $X$ is said to be equal to $k$ \cite[Definition 6.2]{M}.
\par Then $\A$ has the {\it cell decomposition property} (CDP) if any $C$-definable subset $X\subset A^n$, $C\subset A$,  can be partitioned into finitely many cells  and for any $C$-definable function $f$ from $X$ to $\A$ there exists a partition of $X$ into finitely many $C$-definable cells $X_{i}$ such that $f\mid_{X_{i}}$ is continuous \cite[Definition 6.3]{M}.
\edfn 
\par In order to state the main result of Mathews, we need to recall what it means for $\A$ to have the {\it local continuity property}. 
\dfn \label{lcp} Assume that $\A$ has {\it finite Skolem functions} \cite[Definition 5.1]{M}. 
Suppose we are given $k$ polynomials $p_{1}[\bar X,Y],\cdots, p_{k}[\bar X,Y]\in A[\bar X,Y]$, which we may suppose to be closed under derivation with respect to $Y$. Let $d$ be the sum of the degrees of these polynomials with respect to $Y$ and set $S_{j}:=\{\bar a\in A: \A\models \exists^{=j} \;y \bigvee_{i=1}^k p_{i}(\bar a,y)=0\}$, $1\leq j\leq d$.
These subsets are definable and since $\A$ has finite Skolem functions, we get that there exist definable functions $F_{j,r}$ with the property:
$$\A\models \forall \bar x\in S_{j}\;\bigwedge_{1\leq r\leq j}(\bigvee_{i}p_{i}(\bar x,F_{jr}(\bar x))=0\,\&\,\bigwedge_{1\leq r\neq s\leq j}F_{j,r}(\bar x)\neq F_{j,r}(\bar x)).$$
Then $\A$ has {\it the local continuity property}  \cite[page 12]{M} if for every $1\leq j\leq d$, for any $\bar x\in Int(S_{j})$ there exists an open ball $\bar x\in B\subset S_{j}$ such that  $F_{j1},\cdots, F_{j,j}$ are continuous on $B$. 
\edfn
\par If $\A$ is a real-closed field or a $p$-adically closed field, then $\A$ has the local continuity property. 
\ft \label{Mat} \cite[Theorem 7.1,\;Lemma 9.8]{M}
Assume that the $\L$-structure $\A$ satisfies the following
\begin{enumerate}
\item $\A$ is a topological system,
\item $Th(\A)$ admits quantifier elimination,
\item $Th(\A)$ has finite Skolem functions,
\item $\A$ satisfies the local continuity property.
\end{enumerate}
Then $\A$ has the CDP.
\par Assume further that $acl^{\L_{A}}$ is a pregeometry on $A$, then $\A$ has
a fibered dimension function and it coincides with the topological dimension. Moreover,  the dimension of a definable subset of $A$ is equal to $1$ iff it is infinite. 
\eft
\par Note that we include in the definition of "fibered" that it is a definable fibered dimension function.
\medskip
\par From now on in the remainder of the paper, we will assume that the $\L$-theory $T_{c}$ has finite Skolem functions and that the models of $T_c$ have the local continuity property. This will allow us to use Mathews's cell decomposition theorem but besides that, 
we will also need the first assumption independently.
\par In particular this will apply to the cases where $T_{c}=RCF,\;RCVF,\;_{p}{CF}$. Note also that these theories satisfy the further hypothesis that $acl^{\La}$ satisfies the exchange and so $acl^{\L_{A}}$ is a pregeometry on $A$.  (Note also that $T_{c}$ has finite Skolem functions  implies that $acl^{\L}=dcl^{\L}$.)
\par In particular, under these hypothesis, given a model $\M$ of $T_{c,D}^*$, we get two dimension functions, both fibered, one denoted by $\L$-dim, on the reduct $\M_{\L}$ described above and another one $t$-dim on $\M$ (see Fact \ref{cor delta dim}). 
\par In case $T_{c}=ACVF_{0}$ the theory of algebraically closed non-trivially valued fields of equicharacteristic $0$, $T_{c}$ has no finite Skolem functions \cite[Lemma 6.6]{HM}, but nevertheless there is a cell decomposition theorem analogous to the above (see for instance \cite[Theorems 4.1, 4.2]{HM}, where it is done in the general setting of $C$-minimal structures). By a former result of A. Robinson, the theory $T_{c}$ admits q.e., and so $acl^{\La}$ is the algebraic closure $ACL$ (in the field sense). Therefore $acl^{\La}$ has the exchange property and it gives rise to a notion of dimension, which coincides with the notion of topological dimension \cite[Proposition 6.3]{HM}. One can also proceed as in \cite{D89} and show that $ACVF_{0}$ is algebraically bounded \cite[Proposition 2.15 and Examples 2.24 (4)]{D89} which implies that there is a fibered dimension function on definable subsets of models of $ACVF_{0}$.
\medskip
\par We begin by the following easy observation that we will use repeatedly.
\lem\label{ge} Let $F$ be a field of characteristic $0$. The formula: $$\bigwedge_{i\in I} p_{i}(\bar x,x_{k})=0$$ is equivalent to a disjunction of formulas of the form: \\either 
\begin{equation}
(p(\bar x,x_{k})=0\;\&\;\partial_{X_{k}} p(\bar x,x_{k})\neq 0\;\&\;\bigwedge_{j\in J} q_{j}(\bar x)=0\;\&\;q(\bar x)\neq 0),
\end{equation}
 or 
 \begin{equation}
 (\bigwedge_{j\in J} q_{j}(\bar x)=0\;\&\;q(\bar x)\neq 0),
 \end{equation}
  where $p_{i}(\bar X,X_{k}), p(\bar X,X_{k})\in F[\bar X,X_{k}]-\{0\}$, $q_{j}(\bar X), q(\bar X)\in F[\bar X]-\{0\}$.
\elem
\pr We prove the Lemma by induction on the sum of the degrees with respect to the variable $X_{k}$ of the polynomials $p_{i}(\bar X,X_{k})$. W.l.o.g. we assume that $\bigwedge_{i\in I} deg_{X_{k}} p_{i}\neq 0$.
\par Let $P:=\{p_{i}(\bar X,X_{k}): \;i\in I\}$. If the set $P$ consists of just one polynomial $p(\bar X,X_{k})$ with non zero degree in the variable $X_{k}$, we take its derivative with respect to $X_{k}.$ We rename the set $P$ by the set obtained by adding to $P$ the polynomial $\frac{\partial p }{\partial X_{k}}(\bar X,X_{k})$. We also replace the formula $p(\bar x,x_{k})=0$ by the disjunction ($p(\bar x,x_{k})=0\;\&\;\frac{\partial p }{\partial X_{k}}(\bar x,x_{k})\neq 0$ or $p(\bar x,x_{k})=0\;\&\;\frac{\partial p }{\partial X_{k}} p(\bar x,x_{k})= 0$). Call this step, the step $0$.
\par Then, we close the set $P$ by applying the generalized Euclidean algorithm (see section 2). Namely given any pair of polynomials $p_{i_{1}},\;p_{i_{2}}\in P$ with $deg_{X_{k}} p_{i_{2}}\leq deg_{X_{k}} p_{i_{1}}$, we perform the Euclidean division: $c(\bar X)^n.p_{i_{1}}(\bar X,X_{k})=p_{i_{2}}(\bar X,X_{k}).q(\bar X,X_{k})+r(\bar X,X_{k})$ for some $n>0$, where $c(\bar X)$ is the leading coefficient of $p_{i_{2}}(\bar X,X_{k})$ viewed as a polynomial in $X_{k}$ and $deg_{X_{k}} r<deg_{X_{k}} p_{i_{2}}$. 
Set $p_{i_{2}}'(\bar X,X_{k}):=p_{i_{2}}(\bar X,X_{k})-c(\bar X).X_{k}^{d}$, where $d:=deg_{X_{k}}(p_{i_{2}})$.
\par We replace the formula $\bigwedge_{i\in I} p_{i}(\bar x,x_{k})=0$ by the following disjunction: $\theta_{1}(\bar x,x_{k})\;\vee\; \theta_{2}(\bar x,x_{k})$, where 
$\theta_{1}(\bar x,x_{k}):=(\bigwedge_{i\in I-\{i_{2}\}} p_{i}(\bar x,x_{k})=0\;\&\;p_{i_{2}}'(\bar x,x_{k})=0\,\&\;c(\bar x)=0)$ and 
$\theta_{2}(\bar x,x_{k}):=(\bigwedge_{i\in I-\{i_{1}\}} p_{i}(\bar x,x_{k})=0\;\&\,c(\bar x)\neq 0\,\&\;r(\bar x,x_{k})=0)$.
\par Corresponding to the formula $\theta_{1}$,  we replace $P$ by the set obtained from $P$ by replacing the element $p_{i_{2}}(\bar X,X_{k})$ by $p_{i_{2}}'(\bar X,X_{k})$ if $deg_{X_{k}} p_{i_{2}}'(\bar X,X_{k})\neq 0$. If not we delete $p_{i_{2}}(\bar X, X_{k})$ from $P$. 
\par Corresponding to the formula $\theta_{2}$, 
if $deg_{X_{k}} r(\bar X,X_{k})\neq 0$, we replace $P$ by the set obtained from $P$ by replacing the element $p_{i_{1}}(\bar X,X_{k})$ by $r(\bar X,X_{k})$ and if $deg_{X_{k}} r(\bar X,X_{k})= 0$ then we replace $P$ by  deleting in $P$ the element  $p_{i_{1}}(\bar X,X_{k})$.
\par Since the sum of the degrees in $X_{k}$ of the polynomials in $P$ decreases, we will obtain eventually either the empty set or a set $P$ consisting of just one polynomial of non zero degree in the variable $X_{k}$. In that last case, we redo the step $0$, where we added a polynomial with strictly smaller degree in $X_{k}$. So we will eventually obtain an empty set $P$ which corresponds to a formula of the form $(2)$.  \qed

\medskip

\prop \label{*dense1} Let $\M\models T_{c,D}^*$ and $K$ a differential subfield of $M$. Given an $\L_{D,K}$-definable set $S$ in $\M$, there exists an $\L_{K}$-definable subset $S^*$ of $S^{alg}$ such that $S$ is included and dense $\pi_{1}(S^*)$. Moreover, there is a finite partition $\Ps_{S^*}$ of  $S^*$ consisting of $\L_{K}$-definable cells $\tilde C$ such that the $\L$-generic tuples of the form $a^{\nabla}$, $a\in S$, are dense in $\pi_{d}(\tilde C)$ with $d=\L$-dim($\tilde C)=\L$-dim$(\pi_{d}(\tilde C))$ and  $\tilde C\subset dcl_{\L_{K}}(\pi_{d}(\tilde C))$.
\eprop
\pr Let $\phi(x)$ be an $\L_{D}$-formula with parameters in $K$ such that $S=\phi(M)$.  
\par Since $T_{c,D}^*$ admits quantifier elimination (Fact \ref{qe}) and because of our assumptions on the relation symbols of $\L$, $\phi(x)$ is equivalent to a finite disjunction $\bigvee_{j\in J} \phi_{j}(x)$, where $\phi_{j}(x):=(\bigwedge_{i\in I_{j}} p_{i}(x)=0\;\&\;\theta_{j}(x))$, $p_{i}(X)\in K\{X
\}$ and $\theta_{j}^*(x^*)$ defines an open subset of some $M^{n_{j}}$, where $n_{j}$ is the length of the tuple $x^*$ (Notation \ref{star}). Let $S_{j}:=\phi_{j}(M)$; we have $S=\bigcup_{j\in J} S_{j}$. From now on, for ease of notation we will replace  $S_{j}^{alg}$ by $S^{alg}$.
\par Using the CDP property of $T_c$ (Fact \ref{Mat}), we are going to describe a procedure replacing $S^{alg}\subset M^{n}$ by a finite union of cells where the differential tuples will be dense and such that any differential tuple that belongs to $S^{alg}$, will belong to one of them. Each of these cells $C$ will have the following description:
\cl \label{des} The differential tuples in $S^{\n}\cap C$ are dense in $C$ and 
\par either $dim(C)=n$ and so $C$ is an open subset of $M^{n}$, 
\par or $dim(C)=0$ and $C$ is a singleton of the form
 $$\bar v:=(u,g_{1}(u),\cdots,g_{n-1}(u)),$$
where $g_{s}(x)$ are rational functions of $x$, $u^{(s)}=g_{s}(u)$, $1\leq s\leq n-1$, and $u, \bar v\in acl(K)$, 
\par or $0<dim(C)=dim(\pi_{n_{1}}(C))=n_{1}<n$ $($depending on $C)$ and $C$ consists of tuples $$\bar v:=(\bar u,f(\bar u),g_{1}(\bar u,f(\bar u)),\cdots,g_{n-n_{1}-1}(\bar u,f(\bar u))),$$ where $\bar u$ belongs to $\pi_{n_{1}}(C)$, $f(\bar y)$ is a definable Skolem function, $g_{s}(\bar y,z)$ rational functions, with $\ell(\bar y)=\ell(\bar u)$, $1\leq s\leq n-n_{1}-1$.
Moreover $f(\bar y)$ is continuous and $g_{s}(\bar y,f(\bar y))$ are well-defined on $\pi_{n_{1}}(C)$. 
For a tuple $\bar u$ of the form $u^{\nabla_{n_{1}}}$, $f(u^{\nabla})=u^{(n_{1}+1)}, g_{s}(u^{\nabla},f(u^{\nabla})))=u^{(n_{1}+s+1)}$, $1\leq s\leq n-n_{1}-1$.
\ecl
\par In case $I=\emptyset$, we have that $S=\theta(M)$ and so $S^{alg}$ is a $\L_{K}$-definable open subset of $M^{n}$ and therefore a cell $C$. 
\par The set of differential points of $S^{\n}$ in $C$ is dense in $C$. (The argument goes as follows \cite[Proposition 3.14]{GP}. Consider the differential polynomial $X^{(n+1)}=0$. Then given any element $\bar a\in C$, the tuple $(\bar a,0)$ satisfies $x_{n+1}=0$. Since $C$ is open, we can find $W\in \V$ such that $\bar a+W\subset C$. Therefore by the axiomatisation of $T_{c,D}^*$, we can find a differential tuple $\bar b\in \bar a+W$, satisfying $b^{(n+1)}=0$. Since $b^{\n_{n}}\in C$, we have that $b\in S$.) In this case we put $S^{alg}$ in $\Ps_{S^*}$.
\par Now assume that $I\neq \emptyset$ and suppose that $S^{alg}\subset M^{n}$.
 We have that $\L$-dim($S^{alg})=m<n$. 
By the CDP property of $T_c$, $S^{alg}$ can be partitioned into a finite union of cells $C_{j}$. To each of these cells $C_{j}$ we will associate a finite union of cells into which the subset of differential tuples $a^{\n}$ with $a\in S$ is dense.
\par $(a)$ First, assume that $\L$-dim$(\pi_{1}(C_{j}))=0$. So, by assumption on the $\L$-dimension function on $\M_{\L}$, $\pi_{1}(C_{j})$ is a singleton, say $\{a\}$, for some $a\in M$ (see Definition 3.1). This set is included in a set definable by a formula 
of the form $\bigvee_{t} \bigwedge_{i} p_{it}(x)=0$, $p_{it}[X]\in K[X]-\{0\}$. In particular $a$ is algebraic over $K$ and so by considering its minimal polynomial, we can express its successive derivatives $D^{(n)}(a)$ as a rational function $g_n(a)$ (with $g_n(X)\in K(X)$) of $a$. Then we check whether $(a,g_1(a),\cdots)$ belongs to $C_{j}$.  
In this case, we keep the $\L_{K}$-formula expressing that ($(a,g_1(a),\cdots)\in C_{ j}$) with $D(a)=g_1(a)$. We put in $\Ps_S^*$ the tuple $(a,g_1(a),\cdots)$, which belongs to $dcl(\{a\})$.
\par $(b)$ Second, assume that $\L$-dim$(\pi_{1}(C_{j}))=1$ and let $k\geq 2$ be minimum such that $\L$-dim($\pi_{k}(C_{j}))<k$. So this implies that
$\L$-dim$(\pi_{k-1}(C_{j}))=k-1.$ 
The subset of $M^{k-1}$ where the dimension of the fibers of $\pi_{k}(C_{j})$ is equal to $0$, is definable and of dimension $k-1$ (\cite[Lemma 9.8]{M}). 
\par Consider the $\L_{K}$-definable subset $\pi_{k}(C_{j})$ of $M^k$.  
It is defined by a formula of the form $(1)$ $\bigvee_{n\in N}(\bigwedge_{i\in I_{n}} p_{i}(\bar x,x_{k})=0\,\&\;\chi_{n}(\bar x,x_{k}))$, where $\ell(\bar x)=k-1$ and $\chi_{n}(\bar x,x_{k})$ defines an open subset of $M^k$.
Note that since $\L$-dim$(\pi_{k}(C_{j}))<k$, for all $n\in N$, $I_{n}\neq \emptyset$, $p_{i}(\bar X,X_{k})\in K[\bar X,X_{k}]\setminus\{0\}$. Moreover, by Lemma \ref{ge}, we may assume that for some $n\in N$, $I_{n}$ is a singleton, say $I_{n}=\{n\}$ and that $\partial_{X_{k}} p_{n}(\bar x,x_{k})\neq 0$ (here we also use that $\L-dim(\pi_{k}(C_{j}))= k-1)$. Moreover since $C_{j}$ is a cell, this holds for all $n\in N$ and so we may assume that the formula $(1)$ has the form $\bigvee_{n\in N} (p_{n}(\bar x,x_{k})=0\,\&\;\partial_{X_{k}} p_{n}(\bar x,x_{k})\neq 0\,\&\,\chi_{n}(\bar x,x_{k}))$.
\par Since $T_{c}$ has finite Skolem functions, we can pick a point in each fiber in a definable manner. Namely, there are finitely many Skolem  functions $f_{s,n}$, $1\leq s \leq d_{n}$, such that $\forall \bar x\forall x_{k}\;(p_{n}(\bar x,x_{k})=0\rightarrow (\bigvee_{s=1}^{d_{n}} p_{n}(\bar x, f_{s,n}(\bar x))=0))$. Now, since $T_{c}$ has CDP, there exists a finite partition $\Ps$ of $\pi_{k-1}(C_{j})$ into definable cells $E$ such that on each element $E$ of the partition $f_{s,n}\restriction E$ is continuous and that for each $n$, there is only such function $f_{s,n}$ such that $\pi_{k}(C_{\ell j})=\bigcup_{E\in \Ps; 1\leq s\leq d_{n}}\{(\bar x, f_{s,n}(\bar x)):\; \bar x\in E\;\&\;p_{n}(\bar x, f_{s,n}(\bar x))=0\}.$ 
\par First, for the elements $E$ of $\Ps$ of dimension $k-1$, we proceed as follows. Let $(\bar a, a_{k})\in \pi_{k}(C_{\ell j})$ with $\bar a\in E$.
So we get $p_{n}(\bar a,a_{k})=0$ and $\partial_{X_{k}} p_{n}(\bar a, a_{k})\neq 0$. Since the cell $E$ has dimension $k-1$, it is open in $M^{k-1}$. 
By the scheme (DL), there is a point $(b,\cdots, b^{(k-1)})$ close to $(\bar a,a_{k})$ (and so we may assume that $(b,\cdots,b^{(k-2)})$ belongs to $E$), such that $p_{n}(b,\cdots,b^{(k-1)})=0$ and $\partial_{X_{k}} p_{n}(b,\cdots,b^{(k-1)})\neq 0$. Moreover 
we get that $f_{s,n}(b,\cdots,b^{(k-2)})=b^{(k-1)}$. Set $b^{\nabla}:=(b,\cdots,b^{(k-2)})$.
\par Since $\partial_{X_{k}} p_{n}(b,\cdots,b^{(k-1)})\neq 0$, we may express $b^{(k)}$ as a rational function $g_{1}$ of $(b,\cdots,b^{(k-1)})$ and so each $b^{(k+i)}$, $0\leq i<n_{j}-1-k$, as rational functions $g_{1+i}$ of $(b,\cdots,b^{(k-1)})$. 
\par Any tuple $b^{\nabla}\in E$ with $\phi(b)$, has the property that the tuple $$(b^{\nabla},f_{s,n}(b^{\nabla}),g_{1}(b^{\nabla},f_{s,n}(b^{\nabla})),\cdots, g_{n_{j}-k}(b^{\nabla},f_{s,n}(b^{\nabla})))\in C_{j}.$$ 
Set $\bar z:=(z_{1},\cdots,z_{k-1})$; we express that:
 $$\bar z\in E\;\;\&\;\;(\bar z,f_{s,n}(\bar z),g_{1}(\bar z,f_{s,n}(\bar z)),\cdots, g_{n_{j}-k}(\bar z,f_{s,n}(\bar z)))\in C_{j}.\;\;\;\;(2)$$
\par The subset of all tuples $(u_{1},\cdots, u_{n_{j}})\in C_{j}$ of that form is a $\L_{K}$-definable subset of $C_{j}$, which can be expressed as a finite disjoint union of cells $\tilde C$.
\par Either, $\pi_{k-1}(\tilde C)$ has dimension $k-1$ and so we have obtained an $\L_{K}$-definable subset of $S^{alg}$ having the property that the differential tuples of the form $b^{\nabla}$ are dense in its $\pi_{k-1}$-projection and such that any differential tuple belonging to that projection has a differential prolongation in $S^{alg}$.
Moreover, $\L-dim(\tilde C)=k-1$ and has the form $(2)$ and thus $\tilde C\subset dcl(\pi_{k-1}(\tilde C))$. We put $\tilde C$ in $\Ps_{S^*}$.
\par Or $\L-dim(\pi_{k-1}(\tilde C))<k-1$ and so we re-apply the above processes $(a)$, $(b)$ to the $\L_{K}$-definable subset $\tilde C$.
In formula $(2)$, we replace the condition $\bar z\in E$ by a finite disjunction of similar conditions (indexed by a finite set $T$) but on tuples of shorter length $\bar z':=(z_{1},\cdots,z_{k_{t}'-1})$, $1\leq k_{t}'<k$, varying in subcells $E'_{t}$, $t\in T$ and $\L-dim(E'_{t})=k'_{t}$:
$$\bar z'\in E_{t}'\;\;\&\;\;(\bar z',f'_{t}(\bar z'),g'_{1}(\bar z',f'_{t}(\bar z')),\cdots, g_{k-k_{t}'-1}(\bar z',f_{t}(\bar z')))\in E.\;\;\;\;(2)'.$$
\par Since the dimension function on non-empty definable subsets is finite, the procedure will terminate and at the end we will obtain the required partition $\Ps_{S^*}$. 
\qed
\medskip
\rem Note that even though in the above construction we could have replaced the set $S^{alg}$ by $S^{alg}\times M$, only the open cells $C$ in the set $S^*$ that we obtained 
would have been replaced by cells of the form $C\times M$. In all the other cases, namely for cells $D$ of $\L$-dimension strictly less than the dimension of the ambient space, the successive derivatives of an element belonging to $D$ belong to $dcl(D)$.
\erem
\medskip
\par Applying our analysis on the relationship between $\L$-definable sets and $\L_{D}$-definable ones, we deduce the
following result. 
\dfn The structure $\M$ satisfies the {\it uniform finiteness} $(UF)$ (or equivalently eliminates the quantifier $\exists^{\infty}$) if for every $m, n\in \IN^*$, for every definable subset $S$ (with parameters) of $M^m\times M^n$, there exists $N\in \IN$ such that either $\{\bar y\in M^n:\; (\bar x, \bar y)\in S\}$ is infinite or is finite of cardinality $\leq N$. (Note that it suffices to demand it for subsets $S\subset M\times M^n$.)
\edfn
\cor \label{UF} Let $\M\models T_{c,D}^*$, then $\M$ satisfies $(UF)$. 
\ecor
\pr Let $S$ be an $\L_{D}$-definable subset of $\M$ (with parameters) and let $S^*$ be the associated $\L$-definable subset (with parameters) in some cartesian product of $M$ such that $S$ is included and dense in $\pi_{1}(S^*)$ (Proposition \ref{*dense1}). 
\par Let $\{S_{\bar a}:\;\bar a\in M^n \}$ be a family of subsets of $M$, defined by the $\L_{D}$-formula $\phi(x,y_{1},\cdots,y_{n})$.  
Let $N$ be the maximal order of the variable $x$ occurring in $\phi$. 
\par By Proposition \ref{*dense1}, the corresponding $\L$-definable set $S_{\bar a}^*$ is a finite union of cells $C_{i,\bar a}$ such that for each $\bar a$, $i\in I(\bar a)$ with $I(\bar a)$ finite. Because the fact that the $\L$-reduct of $\M$ satisfies the CDP, is a property of the $\L$-theory of $\M_{\L}$, we get that the number of cells is bounded independently of $\bar a$, namely for some constant $f$, we have $\vert I(\bar a)\vert\leq f$. 
\par Now, assume that for some $\bar a$, $S_{\bar a}$ is finite, so the corresponding set $S_{\bar a}^*$ is a union of $\leq f$ cells $C_{i,\bar a}$ of $\L$-dimension $0$. Therefore such cell $C_{i,\bar a}$ consists of a point and so we get the bound $f$ on the cardinality of finite fibers $S_{\bar a}$.
\qed
\medskip
\par It is known that a dp-minimal field eliminates $\exists^{\infty}$ \cite[Observation 9.2.10]{JW}, but note that, for instance $CODF$ is not dp-small \cite[Chapter]{Br}. 
\cor {\rm \cite[Corollary 9.2.11]{JW}} Let $\M\models T_{c,D}^*$, then any infinite externally definable subset of $M$ contains an infinite internally definable subset of $M$.
\ecor
\pr It follows from the the fact $T_{c,D}^*$ is NIP and eliminates $\exists^{\infty}$. \qed

\medskip
\par In the following proposition we will use the notations $S$ and $S^*\subset S^{alg}$ introduced in Proposition \ref{*dense1}.
\prop \label{*function1} Let $\M\models T_{c,D}^*$ and $K$ a differential subfield of $M$. To any $\L_{D,K}$-definable unary function $f$ with domain $S\subset M$, we can associate a $\L_{K}$-definable function $f^*$ defined on a large subset of $S^*\subset S^{alg}$ such that  
 we have for all $u\in S$ that $f(u)=f^*(u^{\n})$. Moreover, for each $\tilde C\in \Ps_{S^*}$ of $\L$-dimension $d$, a large subset of the projection $\pi_{d}(\tilde C)$ can be partitioned in a finite union of  $\L_{K}$-definable cells $C$ such that $f^*\restriction {C}$ is continuous.
\eprop
\pr Let $\phi(x,y)$ be an $\L_{D,K}$-formula such that $f(x)=y$ if and only if $\M\models \phi(x,y)$.  
\par Let $S:=dom(f)\subset M$, namely $S=\{x\in M:\exists y\;\phi(x,y)\}$ and assume that $S^{alg}\subset M^s$. 
In Proposition \ref{*dense1}, we showed how to associate with $S$ an $\L$-definable subset $S^*$ and a finite partition $\Ps_{S^*}$ with the following properties: $S^*$ is a finite (disjoint) union of $\L_{K}$-definable cells $\tilde C\subset M^{n(\tilde C)}$, $n(\tilde C)\leq s$, $\tilde C\in \Ps_{S^*}$ such that, assuming that $\L$-dim$(\tilde C)=d(\tilde C)$, and setting $\bar x_{d(\tilde C)}:=\pi_{d(\tilde C)}(\bar x)$, then, with for ease of notation $d:=d(\tilde C)$,
\par  either $\{b^{\nabla}: b^{\nabla}\in \tilde C\}$ is dense in $\tilde C$, ($d=n(\tilde C)$); we set $\phi_{\tilde C}^*(\bar x_d,\bar y):= 
(\bar x_d\in \tilde C\;\&\,\phi^*(\bar x_d,\bar y))$, \par  or 
 $\{b^{\nabla}: (b^{\nabla},g(b^{\nabla}),g_{1}(b^{\nabla},g(b^{\nabla})),\cdots, g_{n(\tilde C)-d-1}(b^{\nabla},g(b^{\nabla})))\in S^{alg}\}$ is dense in $\pi_{d}(\tilde C)$, $d<n(\tilde C)$, $g$ is a definable Skolem function continuous on $\pi_{d}(\tilde C)$ and $g_{1},\cdots,g_{n(\tilde C)-d-1}$ are rational functions defined over $\pi_{d+1}(\tilde C)$;
 we set \\$\phi_{\tilde C}^*(\bar x_d,\bar y):=
(\phi^*(\bar x_d, g(\bar x_d),g_{1}(\bar x_{d},g(\bar x_{d})),\cdots,g_{n(\tilde C)-d-1}(\bar x_{d},g(\bar x_{d})),\bar y)\;\&\,\bar x_{d}\in \pi_{d}(\tilde C)).$ Moreover $S$ is dense in $\bigcup_{\tilde C\in \Ps(S^*)} \pi_{1}(\tilde C)$.
Let $\phi^*_{mod}(\bar x,\bar y):=\bigvee_{\tilde C\in \Ps_{S^*}}\phi^*_{\tilde C}(\bar x_{d(\tilde C)},\bar y)$.
\par By the quantifier elimination result and by Lemma \ref{ge}, each $\phi^*_{\tilde C}(\bar z,\bar y)$ is equivalent to a finite disjunction of the form:
$\bigvee_{n\in N} \phi^*_{\tilde C,n}(\bar z,\bar y)$, where $\phi^*_{\tilde C,n}(\bar z,\bar y):=$
$$p_{n}(\bar z,y_{1},\cdots,y_{m_{n}})=0\;\&\;\partial_{X_{d+m_{n}}} p_{n}(\bar z,y_{1},\cdots,y_{m_{n}})\neq 0\;\&\;\theta_{n}(\bar z,y_{1},\cdots,y_{m_{n}}),$$ where $\theta_{n}^*(M)$ is an open subset and $N$ a finite subset of $\IN$. We allow the possibility for $m_{n}$ to be zero if no $y_{i}$, $1\leq i\leq m_{n}$ occur non-trivially. 
\cl For each $\tilde C\in \Ps_{S^*}$, (setting $d:=d(\tilde C)$) for any $\bar x_d$ on a large subset of $\tilde C$, we have: $$\L-dim(\pi_{1}(\{\bar y:\;\phi^*_{\tilde C}(\bar x_d,\bar y)\}))=0.$$
\ecl
\prcl Indeed, the elements of the form $x^{\n_{d}}$ are dense in $\pi_{d}(\tilde C)$ and for such tuple there exists a unique $y$ such that $\phi_{mod
}^*(x^{\n},y^{\n})$ holds $(\star)^1$. Suppose otherwise and let $u\in M$ be such that for some $n\in N$, $$\L-dim(\pi_{1}(\{\bar y: \phi_{\tilde C,n}^*(u^{\n},\bar y)\})>0.$$ Then, $\pi_{1}(\{\bar y: \phi_{\tilde C,n}^*(u^{\n},\bar y)\})$ would contain an open subset of $M$ with two distinct elements $a_{1}\neq b_{1}$ such that for some $a_{2},\cdots, a_{m_{n}},\;b_{2},\cdots,b_{m_{n}}$ we have\\
$\phi_{\tilde C,n}^*(u^{\n},a_{1},a_2,\cdots,a_{m_{n}})\,\&\,\phi_{\tilde C,n}^*(u^{\n},b_{1},b_2,\cdots,b_{m_{n}})$, namely 
\par $p_{n}(u^{\n},a_{1},\cdots,a_{m_{n}})=0\;\&\;\partial_{X_{n(\tilde C)+m_{n}}} p_{n}(u^{\n},a_{1},\cdots,a_{m_{n}})\neq 0\;\&\;\theta_{n}(u^{\n},a_{1},\cdots,a_{m_{n}})$ and $p_{n}(u^{\n},b_{1},\cdots,b_{m_{n}})=0\;\&\;\partial_{X_{n(\tilde C)+m_{n}}} p_{n}(u^{\n},b_{1},\cdots,b_{m_{n}})\neq 0\;\&\;\theta_{n}(u^{\n},b_{1},\cdots,b_{m_{n}})$.
Applying the scheme $(DL)$, we can find two differential tuples $a^{\n},\;b^{\n}$ close to $(a_{1},\cdots,a_{m_{n}})$, respectively to 
$(b_{1},\cdots, b_{m_{n}})$) such that $\phi_{\tilde C,n}^*(u^{\n},a^{\n})$ and $\phi_{\tilde C,n}^*(u^{\n},b^{\n})$ hold, contradicting $(\star)^1$. \qed

\par Then, we consider the formula $\psi_{\tilde C,f}(\bar x_{d})$ expressing that $\L-dim(\pi_{1}(\{\bar y: \phi^*_{\tilde C}(\bar x_{d},\bar y)\}))=0.$ 
As shown above, the complement of this set in $\pi_{d}(\tilde C)$ is of dimension $<d$. (Otherwise, we could find in it an element of the form $u^{\n}$ and obtain a contradiction by the above.) The same argument shows also that $S$ is dense in $\bigcup_{\tilde C\in \Ps(S^*)} \pi_{1}(\psi_{\tilde C,f}(M))$.
\par Let us denote by $\chi_{\tilde C}(\bar x_{d},y_{1})$ a formula expressing that $\psi_{\tilde C,f}(\bar x_{d})\;\&\; \exists y_{2}\cdots\exists y_{s}\phi^*_{\tilde C}(\bar x_{d},\bar y).$
Again by q.e. and Lemma \ref{ge}, the formula $\chi_{\tilde C}(\bar z,y)$ is of the form $\bigvee_{i\in I} \chi_{\tilde C,i}(\bar z,y)$ with $\chi_{\tilde C,i}(\bar z,y):=(q_{i}(\bar z,y)=0\;\&\;\partial_{X_{n(\tilde C)+1}} q_{i}(\bar z,y)\neq 0\;\&\;\rho_{i}(\bar z,y))$, where $\rho_{i}^*(M)$ is an open subset and for all $i\in I$, $q_{i}(\bar Z,Y)$ is a non-trivial polynomial.  
\par Again we use that $T_{c}$ has finite Skolem functions. So we obtain a finite set $\Fe_{\tilde C}$ of Skolem functions $h$ such that for any tuple of elements $(\bar a,b)$ such that $\chi_{\tilde C}(\bar a,b)$ holds, we can express $b$ as $h(\bar a)$.
Moreover using that for some $i\in I$, $\chi_{\tilde C,i}(\bar a,b)$ holds, we express the successive derivatives of $b$ as rational functions of $(\bar a,b)$. 
Then we express that the resulting tuple satisfies $\phi^*_{\tilde C}$, namely that $\xi_{\tilde C,h}(\bar z):=\phi^*_{\tilde C}(\bar z,h(\bar z),\cdots)$ holds. 
\par For a tuple of the form $a^{\n}$, we have only one function $h\in \Fe_{\tilde C}$ such that $\xi_{\tilde C,h}(a^{\n})$ holds since it implies that $a\in dom(f)$ $(\star)^2$.
 For each of these functions $h\in \Fe_{\tilde C}$, we consider the subset $S_{h}$ of $\pi_{d}(\tilde C)$ such that the formula $\xi_{\tilde C,h}(\bar z)$ holds. Note that if $h\neq h'$, $\L-dim(S_{h}\cap S_{h'})<d$, otherwise we contradict $(\star)^2$. 
\par  Moreover by the same reasoning, $\bigcup_{h\in \Fe_{\tilde C}}\pi_{1}(\xi_{\tilde C,h}(M))$ is dense in $\pi_{1}(\tilde C)$ and if $u\in \pi_{1}(\xi_{\tilde C,h}(M))$, we get that $h(u^{\n})=f(u)$.
\par We get a partition of a large subset of $\pi_{d}(\tilde C)$ that we refine further, using the CDP theorem, into cells on which each $h\in \Fe_{\tilde C}$ is continuous.
\par So for $x^{\n_{d}}$ in each of these sub-cells of $\pi_{d}(\tilde C)$, we get that $\bigvee_{h\in \Fe_{\tilde C}} h(x^{\n_{d}})=f(x)$. 
 \qed
\medskip
\rem Again one can ask what would happen if we associate to the domain $S$ of the function $f$, the set $S^{alg}\times M$ instead of $S^{alg}$. Recall that $\phi(x,y)$ denoted the formula defining the graph of $f$ and that we associated the a cell decomposition of $S^{alg}$ $\L$-formulas of the form $\phi^*_{\tilde C}(\bar x_{d(\tilde C)},\bar y)$. Moreover each of them was equivalent to a finite disjunction of formulas of the form $p_n(\bar x_{d(\tilde C)},\bar y)=0\;\&\;...$. If we apply $D$ to 
$p_n(x^{\n_d},y^{\n_d'})$ we obtain a polynomial $q_n(x^{\n_{d+1}},y^{\n_{d'+1}})$. So we may add the condition $q_n(\bar x_{d_{\tilde C}},x_{d_{\tilde C}+1},\bar y, y_{n_m+1})=0$. As before either belonging to the cell $C$ induces a condition on the successive derivatives of $x$, or $C$ is an open cell.
\erem
\medskip
\prop \label{*densek} Let $\M\models T_{c,D}^*$ and $K$ a differential subfield of $M$. Given an $\L_{D,K}$-definable set $S\subset M^k$, $k\geq 2$, there exists an $\L_{K}$-definable subset $S^*\subset M^{n_{1}+\cdots+n_{k}}$ of $S^{alg}$ such that $S^{\n}$ is dense in $S^*$.
Moreover, $S^*$ is a finite union of $\L_{K}$-cells $\tilde C$ with the following property. 
The projection $\pi_{\bar c_{11}\bar 0\cdots\bar 0\bar c_{k1}}(\tilde C)$ is an open subset of $M^{n_{11}+\cdots+n_{k1}}$, with $\bar c_{\ell1}$, $1\leq \ell\leq k$, a tuple of $1$'s of length $n_{\ell1}\leq n_{\ell}$ 
and each $\tilde C$ is included in $\L_{K}-dcl(\pi_{\bar c_{11}\bar 0\cdots\bar 0\bar c_{k1}}(\tilde C))$.
\eprop
\pr Let $\phi(u_{1},\cdots,u_{k})$ be an $\L_{D,K}$-formula such that $S=\phi(M)$, $k\geq 1$.  
\par Since $T_{c,D}^*$ admits quantifier elimination (see Fact \ref{qe}), $\phi(u_{1},\cdots, u_{k})$ is equivalent to a finite disjunction $\bigvee_{j} \phi_{j}(u_{1},\cdots,u_{k})$, where $$\phi_{j}(u_{1},\cdots,u_{k}):=\bigwedge_{i\in I_{j}} p_{i}(u_{1},\cdots,u_{k})=0\;\&\;\theta_{j}(u_{1},\cdots,u_{k}),$$ $p_{i}(X_{1},\cdots,X_{k})\in K\{X_{1},\cdots,X_{k}\}$ and $\theta_{j}^*(\bar u_{1},\cdots, \bar u_{k})$ defines an open subset of $M^{N_{j}}$, where $n_{\ell}:=\ell(\bar u_{\ell})$, $1\leq \ell\leq k$ and $N_{j}:=\sum_{\ell=1}^k n_{\ell}$. Let $S_{j}:=\phi_{j}(M)$, $S_{j}^{alg}:=\phi_{j}^*(M)$. 
\par For the indices $j$ such that $I_{j}=\emptyset$, we get $S_{j}^{*}=\theta_{j}^*(M)=S_{j}^{alg}$ and $dim(S_{j}^{*})=N_{j}$. By the density property of differential tuples and the fact that $S_{j}^{*}$ is open, we have that $S^{\n}\cap S_{j}^*$ is dense in $S_{j}^*$.

\par Now assume that $j$ is such that $I_{j}\neq \emptyset$ and assume that $\L$-dim($S_{j}^{alg})=d_{j}<N_{j}$. From now on we will drop the index $j$. By the CDP theorem for models of $T_{c}$ (Fact \ref{Mat}), $S^{alg}$ is a finite disjoint union of cells $C\subset M^{N}$ in the $\L$-reduct of $\M$. 
\par For each cell $C$, there is a projection map which is a definable homeomorphism to an open subset. We choose such projection maps $\pi_{\bar c_{1}\cdots\bar c_{k}}$, with $\bar c_{s}$, $1\leq s\leq k$, consisting of $0$ and $1$'s, as follows: write $\bar c_{s}:=\bar c_{s1}\bar c_{s2}$, $1\leq s\leq k$, where $\bar c_{s1}$ consists only of $1$'s (with length possibly $0$) and $\bar c_{s2}$ begins with $0$, then we choose $\pi_{\bar c_{1}\cdots\bar c_{k}}$ such that the tuple $\bar c_{11}0\bar c_{21}0\cdots\bar c_{k1}0$ is maximum in the lexicographic ordering $(\star)$. Let $\ell(\bar c_{s1})=n_{s1}$. As in Proposition \ref{*dense1}, we will define a subset $C^*$ of $C$ such that $C\cap S^{\n}$ is dense in $C^*$.
We proceed by induction on $k$. 
\par We use the description obtained in Claim \ref{des} in Proposition \ref{*dense1} to $\pi_{n_{1}}(C)$. (Note that we can express $\pi_{n_{1}}(C)$ as $S_{1}^{alg}$ for some $\L_{D}$-definable set $S_{1}$.) So we can express this projection as a finite (disjoint) union of subsets 
which can be described as follows. 
\par Either, $\ell(\bar c_{11}):=n_{11}$ is equal to $n_1$, and so $\pi_{n_1}(C)$ is an open subset of $M^{n_{1}}$, or 
  $0\leq n_{11}< n_{1}$, and $\pi_{n_1}(C)$ consists of the tuples $\bar v_{1}\in \pi_{n_{1}}(C)$ 
  \par either of the form, in case $\ell(\bar c_{11})=0$,
$$\bar v_{1}:=(u_{1},g_{1}(u_{1}),\cdots,g_{n_{1}-1}(u_{1})),$$
where $g_{s}(u_{1})$ are rational functions and are equal to the successive derivatives of $u_{1}$, $1\leq s\leq n_{1}-1$, $u_{1}\in acl(K)$, 
\par or of the form, in case $\ell(\bar c_{11})>0$,  $$\bar v_{1}:=(\bar u_{1},f_{1}(\bar u_{1}),g_{11}(\bar u_{1},f_{1}(\bar u_{1})),\cdots,g_{1n_{1}-n_{11}}(\bar u_{1},f_{\ell}(\bar u_{1}))),$$ where $f_{1}$ is a definable Skolem function of $\bar u_{1}$, $g_{1s}(\bar u_{1})$ rational functions, $1\leq s\leq n_{1}-n_{11}-1$, $\bar u_{1}:=(u_{11},\cdots, u_{1n_{11}})$ belongs to a finite union of subcells $E$ of $\pi_{n_{11}}(C)$.
Moreover $f_{1}$ is continuous on each subcell $E$ and for a tuple $\bar u_{1}$ of the form $u_{1}^{\nabla}$, $f_{1}(u_{1}^{\nabla}), g_{1s}(u_{1}^{\nabla}, f_{1}(u_{1}^{\nabla}))$ are its successive derivatives.
\par To ease notation we will abbreviate such tuple by $(\bar u_1, F_1(\bar u_1))$, where $F_1$ denotes the tuple of $\L_{K}$-definable functions occurring above (namely either a Skolem function or the composition of a rational function with a Skolem function with coefficients in $K$) and to treat the first case (when $n_{11}=n_1$) in the same way, we allow the tuple of functions $F$ to be empty.
\cl\label{desk} To each cell $C$, by induction on $1\leq s\leq k$, we associate with $\pi_{\bar c_{1}\cdots\bar c_{s}}(C)$, an $\L_{K}$-definable set $C_{s}^*$ which is itself a finite union of $\L_{K}$-definable subsets of the form $C_{1,2,\cdots,s}$, consisting of tuples of the form: $$(\bar u_{1},F_{1}(\bar u_{1}),\bar u_{2},F_2(\bar u_1,\bar u_2),\cdots,\bar u_{s}, F_{s}(\bar u_{1},\bar u_{2},\cdots, \bar u_{s})),$$
 where $(\bar u_{1},\cdots,\bar u_{s})$ vary in open subcells $\tilde C$ of
$\pi_{\bar e_{11}\cdots{\,}\bar e_{s1}}(C)$, for some tuples $\bar e_{11},\cdots, \bar e_{s1}$ of $1$'s, with $\bar c_{11}={\bar e_{11}}^{\smallfrown}\bar 1,\cdots, \bar c_{s1}=\bar e_{s1}^{\smallfrown}\bar 1$ with the possibility that $\ell(\bar e_{i1})=0$, $1\leq i\leq s$, in which case the corresponding $u_{i}$ belongs to $acl(K,\bar u_{1},\cdots,\bar u_{i-1})$ (with the convention that $acl(K,\bar u_{0})=acl(K)$),
$F_i(\bar u_1,\cdots,\bar u_i)$ are tuples of Skolem or composition of rational functions and Skolem functions as above, $1\leq i\leq s$, with the convention that if $n_{i1}=n_i$, the tuple $F_i$ is empty. Moreover the $F_{i}$ are well-defined and continuous on each $\tilde C$ and for tuples $\bar u_{i}$ of the form $u_{i}^{\n}$, $1\leq i\leq s$, the tuple $(u_{i}^{\n},F_{i}(u_{1}^{\n},\cdots,u_{i}^{\n}))$ is equal to the successive derivatives of $u_{i}$ up to order $n_{1}+\cdots+n_{i}-1$.
\ecl
\pr By the discussion above, the Claim holds for $s=1$. Assume it holds for indices up to $1\leq s<k$ and consider $\pi_{\bar c_{1}\cdots\bar c_{s+1}}(C)$. We fix an element $\bar v_{s}\in E$, where $E$ is one of the subcells $C_{1,\cdots,s}$ of $\pi_{\bar c_{1}\cdots\bar c_{s}}(C)$ occurring in Claim \ref{desk}, treating that element as a parameter. 
\par Either, $\ell(\bar c_{s+1\;1})=n_{s+1}$ and so $\pi_{\bar 0\bar c_{s+1}}(C)$, with $\ell(\bar 0)=n_{1}+\cdots+n_{s}$, is an open subset of $M^{n_{s+1}}$.  We define $C_{1,\cdots, s+1}$ as $\{(\bar v_{s},\bar u_{s+1})\in\pi_{n_{1}+\cdots+n_{s+1}}(C):\;\bar v_{s}\in C_{1,\cdots,s}\}$.
\par Or, 
 $n_{s+1\;1}<n_{s+1}$; we express that $(\bar v_{s},u_{s+1})\in\pi_{n_{1}+\cdots+n_{s}+1}(C)$ for $\bar v_{s}\in E$ and \\
 in case $n_{s+1\;1}=0$, we obtain tuples of the form $$(\bar v_{s},u_{s+1},g_{s+11}(\bar v_{s},u_{s+1}),\cdots,g_{s+1\;n_{s+1}-1}(\bar v_{s},u_{s+1})),$$ 
with $u_{s+1}\in acl_{\L}(K,\bar v_{s})$, plus possibly an additional open condition on $(\bar v_{s},u_{s+1})$ of the form $q(\bar v_{s},u_{s+1})\neq 0$ with $q[\bar X]\in K[\bar X]\setminus \{0\}$; \\
in case $0<n_{s+1\;1}$, we obtain tuples of the form
 $$(\bar v_{s},\bar u_{s+1}, f_{s+1}(\bar v_{s},\bar u_{s+1}),g_{s+1\;1}(\bar v_{s},\bar u_{s+1},f_{s+1}(\bar v_{s},\bar u_{s+1})),\cdots,g_{s+1\;n_{s+1}-n_{s+1\;1}-1}(\bar v_{s},\bar u_{s+1},f_{s+1}(\bar v_{s},\bar u_{s+1}))),$$ where $\bar u_{s+1}:=(u_{s+1\;1},\cdots, u_{s+1\;n_{s+1\;1}})$, $f_{s+1}$ is a definable Skolem function, $g_{s+1\;i}(\bar v_{s},\bar u_{s+1},x)$ rational functions, $1\leq i\leq n_{s+1}-n_{s+1\;1}-1$ and $\bar u_{s}$ belongs to $\pi_{\bar \bar 0\bar c_{s1}}(C)$ with $\ell(\bar 0)=n_{1}+\cdots+n_{s}$. Moreover we may refine $E$ in such a way that $f_{s+1}$ is continuous on $E\times \tilde E$, with $\tilde E$ an open subcell of $\pi_{\bar \bar 0\bar c_{s1}}(C)$, and the rational functions $g_{s}$ well-defined.
 \medskip
\par Finally we express that the tuples (obtained in Claim \ref{desk}) $(\bar u_{1},F_{1}(\bar u_{1}),\cdots,\bar u_{k}, F_{k}(\bar u_{1},\bar u_{2},\cdots, \bar u_{k}))$ belongs to $C$ which may impose further conditions on $(\bar u_{1},\cdots,\bar u_{k})$. On subcells of the same $\L$-dimension as $C$, we can stop. On  subcells of smaller dimension than $\L$-dim$(C)$, we re-iterate the procedure (which will eventually terminate). 
\qed
\medskip
\rem Note that the $\L$-definable set $S^*$ that we associated with an $\L_{D}$-definable set $S$ has the property that if $S^{alg}$ has a non empty interior, then $S^*$ has a non-empty interior.
\erem
\medskip
\par Keeping the notation of the previous proposition, and in particular denoting by $\Ps_{S}^*$ the set of all cells of the form $\tilde C$ with the properties stated in Proposition \ref{*densek}, we get:
\prop \label{*functionk} Let $\M\models T_{c,D}^*$ and $K$ a differential subfield of $M$. To any $\L_{D,K}$-definable function $f:M^k\rightarrow M^k$ with domain $S\subset M^k$, we can associate an $\L_{K}$-definable function $f^*$ defined on a large subset of $S^*\subset S^{alg}$ such that for any $u\in S$, we have $f(u)=f^*(u^{\n})$. Moreover, a large subset of each $\pi_{\bar c_{11}\bar 0\cdots\bar c_{k1}}(\tilde C)$, with $\tilde C\in \Ps_{S^*}$, can be partitioned in a finite union of  $\L_{K}$-definable cells $C$ such that $f^*\restriction {C}$ is continuous.
\eprop
\pr It suffices to consider the case where the image of $f$ is included in $M$ (we can express $f$ as a $k$-tuple of $\L_{D,K}$-definable functions with values in $M$). Let $\phi(\bar x,y)$ be an $\L_{D,K}$-formula such that $f(\bar x)=y$ if and only if $\M\models \phi(\bar x, y)$.  
\par Let $S:=dom(f)$, equivalently $\{\bar x\in M^k:\exists y\;\phi(\bar x,y)\}$. In Proposition \ref{*densek}, we associated with $S$ a $\L_{K}$-definable subset $S^*$ and a finite partition $\Ps_{S^*}$ consisting of $\L_{K}$-definable sets $\tilde C$ of the form described in Claim \ref{desk}: it consists of tuples $$(\bar u_{1},F_{1}(\bar u_{1}),\bar u_{2},F_2(\bar u_1,\bar u_2),\cdots,\bar u_{k}, F_{k}(\bar u_{1},\bar u_{2},\cdots, \bar u_{k})),$$
 where $\bar u_{i}:=(u_{i1},\cdots, u_{in_{i1}})$, $0\leq n_{i1}\leq n_{i}$, with $F_i(\bar u_1,\cdots,\bar u_i)$ tuples of definable continuous Skolem or rational functions following the same convention as above, $1\leq i\leq k$ and $(\bar u_{1},\cdots,\bar u_{k})$ varies in the open set
$\pi_{\bar e_{11}\bar 0\cdots{\,}\bar e_{k1}}(\tilde C)$, where $\bar e_{i1}$ is a tuple of $1$'s of length $n_{i1}$. For ease of notation let us denote $\pi_{\bar e_{11}\bar 0\cdots{\,}\bar e_{k1}}$ by $\pi$.
\par Set $\psi_{\tilde C}(\bar u,\bar y):=(\phi^*((\bar u, F(\bar u)),\bar y)\;\&\;(\bar u,F(\bar u))\in \tilde C),$ with $\ell(\bar u)=\sum_{i=1}^k n_{i1}:=n_{\tilde C}$. W.l.o.g. we may assume that $\L$-dim$(\tilde C)>0$.
\par The proof is along the same lines as the proof of Proposition \ref{*function1}, replacing the tuple $\bar u$ by the tuple $(\bar u_{1},\cdots, \bar u_{k})$. For each $\tilde C$, one shows that on a large subset of $\tilde C$, we may define a function which coincide with $f$ on the subset of differential tuples. Moreover we may definably partition this large subset such that on each element of this (finite) partition, there will be a continuous $\L_{K}$-definable function equal to $f$ on the subset of differential tuples.
\cl\label{large} The set $\bar u\in \pi(\tilde C)$ for which $\exists \bar y\;\psi_{\tilde C}(\bar u,\bar y)$ holds, forms a large subset of $\pi(\tilde C)$. %
\ecl
\pr We first show that the set $\bar u\in \pi(\tilde C)$ for which $\exists\bar y\;\psi_{\tilde C}(\bar u,\bar y)$ holds, forms a large subset of $\pi(\tilde C)$.
Otherwise, we could find in the complement 
a tuple of the form $\bar u^{\n}$ such that $\forall \bar y\;\neg\;\phi^*(\bar u^{\n},\bar y)$, which would contradict that $\tilde C\in \Ps_{S}^*$. 
 \qed
\medskip
\par Let $\psi_{\tilde C,1}(\bar u,\bar y):=\psi_{\tilde C}(\bar u,\bar y)\;\&\;\exists\bar y\;\psi_{\tilde C}(\bar u,\bar y)$. 
\cl\label{funct} For each $\tilde C\in \Ps_{S}^*$, let us show that for any $\bar u$ on a large subset of $\pi(\tilde C)$, \\$\L-dim(\pi_{1}(\{\bar y:\;\psi_{\tilde C,1}(\bar u,\bar y)\}))=0$.
\ecl
\pr Indeed the elements of the form $\bar u^{\n}$ are dense in $\pi(\tilde C)$ and $\exists^{=1} y\;\phi^*(\bar u^{\n}, y^{\n})$ holds $(1)$.
\par By the fact that $T_{c}$ admits quantifier elimination and by Lemma \ref{ge}, $\psi_{\tilde C, 1}(\bar u,\bar y)$ is equivalent to a finite disjunction of the form: $\bigvee_{n\in N} \varphi_{n}(\bar u,\bar y)$, where $\varphi_{n}(\bar u,\bar y)$ is either of the form:
$$(i)\;\;p_{n}(\bar u,y_{1},\cdots,y_{m_{n}})=0\;\&\;\partial_{X_{n_{\tilde C}+m_{n}}} p_{n}(\bar u,y_{1},\cdots,y_{m_{n}})\neq 0\;\&\;\theta_{n}(\bar u,y_{1},\cdots,y_{m_{n}}),$$ or of the form: 
$$(ii)\;\;\theta_{n}(\bar u,y_{1},\cdots,y_{m_{n}}),$$
where $\theta_{n}^*(M)$ is an open subset. We allow the possibility for $m_{n}$ to be zero if no $y_{i}$, $1\leq i\leq m_{n}$ occur non-trivially. By Claim \ref{large}, note that for some $n\in N$, $m_{n}\neq 0$.
Suppose that for some $\bar u$, $\L$-dim$(\pi_{1}(\{\bar y: \psi_{\tilde C,1}(\bar u^{\n},\bar y)\}))>0.$ Then, $\pi_{1}(\{\bar y:  \psi_{\tilde C,1}(\bar u^{\n},\bar y)\})$ would contain an open subcell of $M$. Let $a_{1}\neq b_{1}$ in that open cell such that for some $n_{1},\,n_{2}\in N$, $\bar a,\;\bar b\in M$ with $\pi_{1}(\bar a)=a_{1}$ and $\ell(\bar a)=m_{n_{1}}$, $\pi_{1}(\bar b)=b_{1}$ and $\ell(\bar b)=m_{n_{2}}$, we have: 
$\varphi_{n_{1}}(\bar u^{\n},\bar a)$ 
and $\varphi_{n_{2}}(\bar u^{\n},\bar b)$.
\par Let us assume that $\varphi_{n}$ is in the form $(i)$ (the case $(ii)$ is straightforward). So, $p_{n_{1}}(\bar u^{\n}, \bar a)=0\;\&\;\partial_{X_{n_{\tilde C}+m_{n_{1}}}} p_{n_{1}}(\bar u^{\n}, \bar a)\neq 0\;\&\;\theta_{n}(\bar u^{\n}, \bar a)$, and similarly for $(\bar u^{\n}, \bar b).$
By the scheme $(DL)$, we can find two differential tuples $a^{\n}$ and $b^{\n}$ close to respectively $\bar a$, $\bar b$ (and so distinct),
such that $\varphi_{n_{1}}(\bar u^{\n},a^{\n})$ and $\varphi_{n_{2}}(\bar u^{\n},b^{\n})$ holds, a contradiction with $(1)$. \qed
\medskip

\par Then we consider the formula $\psi_{\tilde C,f}(\bar u)$ expressing that $\L-dim(\pi_{1}(\{\bar y: \psi_{\tilde C,1}(\bar u,\bar y)\})=0.$ 
(This is expressible in a first-order way since the $\L$-dimension coincides with the topological dimension and we have a definable basis of neighbourhoods of $0$.)
\par Let us write the formula $\psi_{\tilde C}(\bar u,\bar y)$ as $\psi_{\tilde C}(\bar u,y,\bar z)$ with $\ell(\bar z)=\ell(\bar y)-1$. Denote by $\chi_{\tilde C}(\bar u,y)$ the formula $\psi_{\tilde C,f}(\bar u)\;\&\; \exists \bar z\;\psi_{\tilde C}(\bar u,y,\bar z).$
Again by q.e. and Lemma \ref{ge}, we may assume that the formula $\chi_{\tilde C}(\bar u,y)$ is of the form $\bigvee_{i\in I} q_{i}(\bar u,y)=0\;\&\;\partial_{X_{n_{\tilde C}+1}} q_{i}(\bar u,y)\neq 0\;\&\;\rho_{i}(\bar u,y)$, where $\rho_{i}^*(M)$ is an open subset and for all $i\in I$, $q_{i}(\bar X,Y)\in K[\bar X,Y]\setminus\{0\}$.
\par Now we use the fact that in models of $T_{c}$, definable sets of dimension $0$ are finite and that we have $\L_{K}$-definable finite Skolem functions. So there is a finite set of definable Skolem functions $\Fe_{\tilde C}$ such that we can express $y$ as $h(\bar u)$, for $h\in \Fe_{\tilde C}$.  Then using the polynomials $q_{i}$, we can express the successive derivatives of $y$ as rational functions $g_{j}$ of $(\bar u,y)$, namely $y^{(j)}=g_{j}(\bar u,y)$, $1\leq j<\ell(\bar y)$. Then we express that for the resulting tuple, $\psi_{\tilde C}(\bar u,h(\bar u),g_{1}(\bar u,h(\bar u)),\cdots,g_{\ell(\bar y)-1}(\bar u,h(\bar u)))$ holds. 
\cl For a large subset of $\tilde C$, we have only one function $h\in \Fe_{\tilde C}$ such that $\chi_{\tilde C}(\bar u,h(\bar u))$ holds.
\ecl
\pr Indeed the statement holds for tuples of the form $\bar u^{\n}$ in $\tilde C$ $(2)$.
For each of these functions $h\in \Fe_{\tilde C}$, we consider the subset $S_{h}$ of $\pi(\tilde C)$ such that the formula $\chi_{\tilde C}(\bar u,h(\bar u))$ holds. Note that if $h\neq h'$, $\L-dim(S_{h}\cap S_{h'})<\L$-dim$(\tilde C)$, otherwise we contradict $(2)$.
\par Using the CDP theorem, we partition each $S_{h}$ into subcells $C_{h}$ on which $h\in\Fe_{\tilde C}$ is continuous. We define on $C_{h}$ an $\L_{K}$-definable function $f_{\tilde C}^*$ equal to $h$. Taking the union over $h\in \Fe_{\tilde C}$, we get a function $f^*$ on $\tilde C$ such that $\pi(dom(f^*))$ is a large subset of $\pi(\tilde C)$.
Moreover for any $\bar u\in S$ such that $\bar u^{\n}\in \tilde C$, we get that $\bigvee_{h\in \Fe_{\tilde C}} (h(\bar u^{\n})=y\rightarrow f(\bar  
u)=y)$. 
 \qed
\medskip

\cor \label{open-def} Let $\M\models T_{c,D}^*$ and $K$ a differential subfield of $M$. Any $\L_{D,K}$-definable open set $S\subset M^k$ is $\L_{K}$-definable subset.
\ecor
\pr By Proposition \ref{*densek}, we may associate with $S\subset M^k$ an $\L_{K}$-definable set $S^*\subset M^{n_{1}+\cdots+n_{k}}$ of $S^{alg}$ such that $S$ is included and dense $\pi_{i_{1},\cdots,i_{k}}(S^*)$, with $i_{1}=1<i_{2}=n_{1}+1\cdots<i_{k}=n_{1}+\cdots+n_{k-1}+1$. So if $S$ is open, then $S=\pi_{i_{1},\cdots,i_{k}}(S^*)$. But $\pi_{i_{1},\cdots,i_{k}}(S^*)$ is $\L_{K}$-definable and so is $S$. \qed 
\medskip
\par Let $\La_{Mac}:=\{+,.,-,0,1,P_{n};\;n\geq 1\}$, let $M$ be a $p$-adically closed field and let $\M$ an $\La$-expansion of $M$. Then $\M$ is $p$-minimal if any $\La$-definable subset of $M$ is already $\La_{Mac}$-definable \cite[Definition 5]{C}. We will make the assumption that when $\M$ is said to be $p$-minimal, then the underlying field is $p$-adically closed.
\cor \label{opencore} Let $\M$ be a model of $T_{c,D}^*$ and suppose its $\L$-reduct is either $o$-minimal or $p$-minimal or weakly $o$-minimal, then its open core is still respectively $o$-minimal, $p$-minimal, weakly $o$-minimal. \qed
\ecor
\medskip

\subsection{Real-closed fields} We will now restrict ourselves to the theory $T_{c}$ of real-closed fields. In particular, in this subsection $\K$ will always denote an expansion of an ordered group endowed with a dense linear order.
\par On one hand, we will give another proof that the open core (Definition \ref{core}) of any model of $CODF$ is o-minimal, using a former result of A. Dolich, C. Miller and C. Steinhorn. Indeed, if $\K$ satisfies the following two properties: uniform finiteness $(UF)$ and definable completeness $(DC)$, then its open core is o-minimal \cite[Theorem A]{DMS}. 
\par Then we will give a proof that $CODF$ admits elimination of imaginaries, which has appeared in \cite{P1}. Recently another proof of elimination of imaginaries has been given in \cite{BCP}, using a fine description of definable types in models of $CODF$ in addition to elimination of imaginaries in $RCF$ and Proposition \ref{*dense1}.
\dfn 
The structure $\K:=(K,<,\cdots)$ is {\it definably complete} $(DC)$ if every unary definable subset (with parameters) of $K$ has both a supremum and an infimum in $K\cup\{\pm\infty\}$.
\edfn 
\par It is well-known that any o-minimal theory expanding a dense order has both properties $(UF)$ and $(DC)$ \cite{D98}. 
\cor {\rm (\cite{P1})}\label{dc} Let $\K$ be a model of $CODF$. Then $\K$ has both properties $(DC)$ and $(UF)$ and so $\K^{o}$ is o-minimal.
\ecor
\pr By Corollary \ref{UF}, $\K$ satisfies $(UF)$, it remains to show that $\K$ satisfies $(DC)$.
\par Let $S$ be an $\L_{D}$-definable subset of $\K$ (with parameters) and let $S^*$ be the associated $\L$-definable subset (with parameters) in some cartesian product of $K$ such that $S$ is included and dense in $\pi_{1}(S^*)$ (Proposition \ref{*dense1}). Since the $\L$-reduct of $\K$ satisfies $(DC)$ and since $\pi_{1}(S^*)$ is $\L$-definable (with parameters), since the supremum and infimum of $S$ are the same as the ones of $\pi_{1}(S^*)$, we get that $\K$ is $(DC)$.
\par This enables us to apply Theorem A in \cite{DMS} and deduce that $\K^{o}$ is o-minimal.
\qed
\medskip
\rem Definably complete fields have been studied by a number of authors and Corollary \ref{dc} allows us to apply their results to any model $\K$ of CODF. Let us mention two of them. For a continuous definable function $g:\K\rightarrow \K$, the intermediate value theorem holds \cite[section 2]{servi}. 
\par Let $f:\K\rightarrow \K$ be a definable monotone function, then $f$ has a derivative $f'$ (with values in $K$), which is continuous except on a finite subset \cite[Theorem B, Lemmas 6.19, 7.1]{FH}. 
\erem
\par Here we will follow the terminology of \cite[Chapter 4.4]{Ho} on (uniform) elimination of imaginaries and also \cite[section 16.5]{Poizat}. Let $T$ be a theory, then $T$ admits elimination of imaginaries (e.i.) if given any model $\M$ of $T$, any formula $\phi(\bar x,\bar a)$ with $\bar a\in M$, there exists a formula $\psi(\bar x,\bar z)$ and a unique tuple $\bar b\in M$ such that $$\M\models \forall \bar x\;(\phi(\bar x,\bar a)\leftrightarrow \psi(\bar x,\bar b))\;\;\;(\star)$$ such tuple $\bar b$ is called the canonical parameter associated to $\phi(M,\bar a)$. The theory $T$ has 1-e.i., whenever the equivalence ($\star$) holds only for tuple $\bar x$ of length $1$.
\par We will apply the fact that $RCF$ eliminates imaginaries \cite[Theorem 4.4.4]{Ho}, where, since there are Skolem functions, one reduces to show that $RCF$ has 1-e.i. \cite[Lemma 4.4.3]{Ho}. Note that in $CODF$, one does not have Skolem functions.
\dfn Given a $(i_{1},i_{2},\cdots,i_{k})$-cell $C$, where $i_{j}\in \{0,1\}$, $0<k$. We will say $C_{1}$ is a {\it subcell} of $C$, if $C_{1}=\pi_{\ell}(C)$ for some $1\leq \ell\leq k$.
\edfn
\par In order to be self-contained, we will give a proof here of the following result. 
\thm {\rm \cite{P1}} The theory $CODF$ eliminates the imaginaries. 
\ethm
 \pr Let $\M\models CODF$ and consider an $\L_{D}$-formula $\phi(x_{1},\cdots,x_{k},\bar a)$, with $\bar a\in M$ and let $S=\phi(M,\bar a)$. 
Let $K$ be the differential subfield of $M$ generated by $\bar a$.
By Proposition \ref{*densek}, there exists a $\L_{K}$-definable subset $S^*\subset M^{n_{1}+\cdots n_{k}}$ of $S^{alg}$ such that $S$ is included and dense $\pi_{i_{1},\cdots,i_{k}}(S^*)$, with $i_{1}=1<i_{2}=n_{1}+1\cdots<i_{k}=n_{1}+\cdots+n_{k-1}+1$.
Moreover, $S^*$ is a finite union of $\L_{K}$-cells $\tilde C$ with the following property $(\star)$: the differential tuples are dense in $\pi_{\bar c_{11}\bar 0\cdots\bar c_{k1}}(\tilde C)$, with $\bar c_{i1}$, a tuple of $1$'s of length $n_{i1}\leq n_{i}$, $1\leq i\leq k$ and $\L$-dim$(\tilde C)=\sum_{i=1}^k  n_{i1}$. Moreover, each $\tilde C$ is included in $\L_{K}-dcl(\pi_{\bar c_{11}\bar 0\cdots\bar c_{k1}}(\tilde C))$. 
\par Since the cells used in Mathews's result (Fact \ref{Mat}) are not necessarily definably connected as the ones used in o-minimal theories \cite[Chapter 3, Proposition 2.9]{D98}, we partition each $\tilde C$ as a finite union of its definably connected components \cite[Chapter 3, Proposition 2.18]{D98}. So we will assume that each $\tilde C$ is maximally definably connected and it will still have property $(\star)$. Since $RCF$ admits the elimination of imaginaries, we may associate to each $\tilde C$ a formula $\xi_{\tilde C}(.,\bar b_{\tilde C})$, where $\bar b_{\tilde C}$ is the canonical parameter associated with $\tilde C$.
\par Then we express by an $\L_{D}$-formula $\chi_{\tilde C}(\bar z)$ \cite[section 6]{KPS} that for some parameters $\bar z$, the definable set $\xi_{\tilde C}(M,\bar z)$ is a (particular) $(\bar k_{1}',\bar k_{2}',\cdots)$-cell with $\bar k_{i}'$ a tuple of length $n_{i}'\leq n_{i}$, of the form $(1,\cdots,1,0,\cdots,0)$ where the differential tuples are dense.
\par Set $\bar b=(\bar b_{\tilde C})$ and let $\xi(\bar x_{1},\bar x_{2},\cdots,\bar x_{k},\bar b):=\bigvee_{\tilde C} (\xi_{\tilde C}(\bar x_{1},\bar x_{2},\cdots,\bar x_{k},\bar b_{\tilde C})\&\;\chi_{\tilde C}(\bar b_{\tilde C}))$.
\par Let $\psi(x_{1},x_{2},\cdots,x_{k},\bar b)$ be the $\L_{D}$-formula expressing that each $\bar x_{j}$, $1\leq j\leq k$, is a differential tuple and that $\xi(\bar x_{1},\bar x_{2},\cdots,\bar x_{k},\bar b)$ holds.
Note that $\psi(M,\bar b)=\phi(M,\bar a)$.
\par Suppose that there exists a tuple $\bar c\neq \bar b$ such that $\psi(M,\bar b)= \psi(M,\bar c).$ Then for some $\bar b_{\tilde C}$ there exists $\bar c_{1}\subset \bar c$ such that $\bar b_{\tilde C}\neq \bar c_{1}$ and $\bar c_{1}$ occurring in $\xi_{\tilde C}(\bar x_{1},\bar x_{2},\cdots,\bar x_{k},\bar c_{1})$. We choose such $\bar b_{\tilde C}$ such that $\tilde C$ is of maximal o-minimal dimension, say $d$.
\par We collect all cells $\tilde C$ of o-minimal dimension $d$ and send each one to one of the cells $D_{1}:=\xi_{\tilde D}(M,\bar c')$, $\bar c'\subset \bar c$ of the same dimension $d$ which has a non-empty intersection. Note that since we have chosen $\tilde C$ to be a definably connected component and since $D_{1}$ is definably connected, it cannot at the same time have a nonempty intersection with $D_{1}$ and its complement, it cannot be strictly included in $D_{1}$ since it is maximally connected and so the only possibility left is that it strictly contains $D_{1}$. Let us show this is impossible.
\par So either the projections by $\pi_{n_{1}'+1}$ of $\tilde C$ and $D_{1}$ differ, or there exists a tuple $\bar x\in \pi_{n_{1}'+1}(\tilde C)$ such that the corresponding fibers are different. Note that if two $(1,\cdots,1,0,\cdots)$-cells differ, they differ by a relatively open set (which cannot be covered by finitely many cells of smaller dimension).
\par In the first case, there exists $\bar x\in \pi_{n_{1}'+1}(\tilde C)\setminus\pi_{n_{1}'+1}(D_{1})$. Since the differential tuples are dense and since $\pi_{n_{1}'+1}(\tilde C)$ is locally closed, there exists a differential tuple $\bar d:=(d,d^{(1)},\cdots,d^{(n_{1}')})$ close to $\bar x$ belonging to $\pi_{n_{1}'+1}(\tilde C)\setminus\pi_{n_{1}'+1}(D_{1})$.
\par In the second case, there is a differential tuple $\bar d$ (in the interior of the projection $\pi_{n_{1}'+1}(\tilde C)$) where the fiber in $\tilde C$ and in $\xi(M,\bar c)$ are different.
\par In both cases, there is a tuple $(\bar d,\bar u)$ in $\tilde C$ and not in $\xi(M,\bar c)$. Again since the differential tuples are dense and since $\tilde C$ is locally closed, there is a differential tuple $(\bar d,\bar e)$ close to $(\bar d,\bar u)$ in $\tilde C$ and not in $\xi(M,\bar c)$, a contradiction. \qed
\medskip
\section{Definable groups}
Throughout this section, we will put the same hypothesis on the theory $T_{c}$ as in sections 2 and 3. In particular, we assume that $T_{c}$ has finite Skolem functions and that the local continuity property holds in models of $T_{c}$.
\par Let $\M$ be a model of $T_{c,D}^*$ and let $(G,.,1)$ be a definable group in $\M$, namely the domain of $G$ is a definable subset of some cartesian product $M^n$ of $M$ and the graph of the group operation $.$ is a definable subset of $M^{3.n}$. We will associate with $G$ a $\L$-definable local subgroup of $\M$. Then we will apply that result to show the descending chain condition on centralisers for definable subgroups in models of CODF. 
\par Let us first recall the notion of a local group \cite[chapter 1, section 20]{MZ}, following the presentation of Goldbring \cite{G}.
\dfn {\rm \cite[section 2]{G}} \label{local-def} A local group $(H,1,i,m)$ is a Hausdorff topological space $H$ endowed with continuous maps $i,\;m$ and an element $1$ such that $i:U\rightarrow H$, where $U$ is an open subset of $H$, $m:O\rightarrow H$ where $O$ is an open subset of $H\times H$ satisfying the following:
\begin{enumerate}
\item $1\in U$, $\{1\}\times H \subset O$, $H\times \{1\}\subset O$,
\item for all $x\in H$, $m(1,x)=m(x,1)=x$
\item for all $x\in U$, we have $(x,i(x))\in O$, $(i(x),x)\in O$ and $m(x,i(x))=1=m(i(x),x)$
\item for all $(x,y),\;(y,z)\in O$ such that $(m(x,y),z)$ and $(x,m(y,z))$ belong to $O$, we have $m(m(x,y),z)=m(x,m(y,z))$.
\end{enumerate}
A {\it local} group is {\it definable} if the above data is given by definable subsets, namely $H,\;U,\;O$ and the graphs of the maps $i,\;m$.
\edfn

\thm \label{local} Let $\M$ be a model of $T_{c,D}^*$, let $(G,\cdot ,^{-1},1)$ be an $\L_{D}$-definable group in $\M$ (possibly with parameters) and let $G^*$ be as in Proposition \ref{*densek}. Then there exists an $\L$-definable (over the same parameters) local group $G^*_{local}$ such that $G^{\n}\cap G^*_{local}$ is dense in $G^*_{local}$ and $\L$-dim$(G^*)$=$\L$-dim$(G_{local}^*)$.
\ethm
\pr Assume that the domain of $G$ is included in some $M^{k}$ and defined by the $\L_{D}$-formula $\phi(u_{1},u_{2},\cdots,u_{k})$, possibly with parameters in $M$, that the graph of the inverse function on $G$ is defined by $\phi_{^{-1}}(u_{1},\cdots,u_{k},v_{1},\cdots,v_{k})$ and the graph of multiplication by $\phi_{\times}(\bar u,\bar v,\bar w)$. We assume that all the parameters occurring in these formulas belong to a subfield $K$ of $M$.
\par In the proof of Proposition \ref{*densek} (Claim \ref{desk}), we showed how to construct from $\phi^*(\bar u^*)$, an $\L$-formula that we denote by $\phi_{mod}^*(\bar u^*)$ by modifying each cell $C$ in $\phi^*(M)$ such that in the new cells, the differential tuples are dense, that each one has only one differential prolongation which belongs to the initial cell and such that any differential tuple belonging to $C$ has a projection in one of these subcells. Assume that the length of the tuple $\bar u^*$ occurring in $\phi^*$ is equal to $m=\sum_{i=1}^k m_{i}$.
\par Let $G^*:=\phi_{mod}^*(M)$ and denote by $1^{\n}$ the differential tuple associated with the neutral element $1$ of $G$. Then in Proposition \ref{*functionk}, we proved that given a definable function and in this case the inverse function on $G$, one can define on a large $\L_{K}$-definable subset $V_{1}$ of $G^*$ a $\L_{K}$-definable function which coincide with $^{-1}$ on $G^{\nabla}$, whose domain contain $G^{\n}$ and which is continuous on a finite definable partition of $V_1$; for ease of notation we will denote this function again by $^{-1}$. 
We keep in that definable partition the cells of the same $\L$-dimension as $V_1$ and we rename the union of these cells $V_1$; note that $V_{1}$ is open but we may be missing some elements of $G^{\n}$ belonging to cells of smaller dimension. However $G^{\n}\cap V_{1}$ is dense in $V_{1}$ and $\L$-dim$(V_{1})=\L$-dim$(G^*)$.
Then we consider the multiplication on $G\times G$ and again using Proposition \ref{*functionk}, we obtain on a large definable subset $Y_{0}$ of $G^*\times G^*$ a $\L_{K}$-definable function that we will denote again by $\cdot$, which coincides with the multiplication of the group $G$ on $G^{\nabla}\times G^{\nabla}$, whose domain contains $G^{\n}\times G^{\n}$ and which is continuous on a finite definable partition of $Y_{0}$. Again, we keep in that definable partition the cells of the same $\L$-dimension as $Y_0$ and we rename the union of these cells $Y_0$ (again this is an open subset with the property that $(G^{\n}\times G^{\n}) \cap O$ is dense and $\L$-dim$(Y)=\L$-dim$(G^*\times G^*)$).
\par  We may assume, by adding that $\phi_{mod}^*(\bar u^*),\phi_{mod}^*(\bar v^*),\phi_{mod}^*(\bar w^*)$ hold, that the definable subsets corresponding to the $\L$-formulas $\phi_{^{-1}}^*(\bar u^*,\bar v^*)$ and $\phi_{\times}^*(\bar u^*,\bar v^*,\bar w^*)$ live in respectively $M^m\times M^m$ and $M^m\times M^m\times M^m$.
\par Consider $V_{0}':=\{\bar y\in V_{1}: \forall \bar a\in G^*$ $\L$-generic of $G^*$ over $\bar y$ $(\bar a,\bar y)\in Y_{0}\;\&\;(1^{\n},\bar y)\in Y_{0}\;\&\;(\bar y,1^{\n})\in Y_{0}\;\&\;\bar y.1^{\n}=1^{\n}.\bar y=\bar y\}$; it is a definable subset of $V_{1}$ by Fact \ref{def-large}. 
\cl  $V_{0}'$ is large in $V_{1}$. 
\ecl
\pr First let us show that $V_{0}:=\{\bar y\in V_{1}: \forall \bar a\in G^*$ $\L$-generic of $G^*$ over $\bar y$ $(\bar a,\bar y)\in Y_{0}\}$ is large and definable in $V_{1}$. Suppose not and let $U_{1}$ be a relatively open subset of $V_{1}\subset G^*$ such that if $\bar y\in U_{1}$, we have that $\{\bar a\in G^*:\;(\bar a,\bar y)\in Y_{0}\}$ is not large. In particular, for all elements of $U_{1}$ of the form $\bar y^{\n}\in V_{1}$ with $\bar y\in G$, we have that the set $\{\bar a:\;(\bar a,\bar y^{\n})\in Y_{0}\}$ is not large. So the complement would contain an open subset $U_{2}$.  But $U_{2}\times U_{1}$ contains tuples of the form $(\bar b^{\n},\bar y^{\n})$ with $\bar b\in G$ which would not belong to $Y_{0}$, a contradiction. 
\par Now let us show that $\{\bar y\in V_{0}: (1^{\n},\bar y)\in Y_{0}\;\&\;(\bar y,1^{\n})\in Y_{0}\}$ is large and definable in $V_{0}$. Suppose not, then there would exist a relatively open subset of $V_{0}$ such that no $\bar y$ in that subset has the property that $(1^{\n},\bar y)\in Y_{0}\;\&\;(\bar y,1^{\n})\in Y_{0}\;\&\;\bar y.1^{\n}=1^{\n}.\bar y=\bar y$. It suffices to take $\bar y$ of the form $\bar z^{\n}$ for some $\bar z\in G$ in order to get a contradiction.
\qed
\medskip
\par Let $V_{0}^{''}:=\{\bar z\in V_{0}':\;\{(\bar x,\bar y)\in Y_{0}:\;(\bar y,\bar z)\in Y_{0}\;\&\;\bar y.\bar z\in V_{0}'\;\&\;(\bar x,\bar y.\bar z)\in Y_{0}\;\&\;(\bar x.\bar y,\bar z)\in Y_{0}\;\&\;\bar x.(\bar y.\bar z)=(\bar x.\bar y).\bar z\}$ is large in $Y_{0}\}$. By the same reasoning as before, the set $V_{0}^{''}$ is definable by Fact \ref{def-large}.
\cl $V_{0}^{''}$ is large in $V_{0}'$. 
\ecl
\pr We proceed by contradiction. If not there would exist $\bar z\in V_{0}'$ and an open set containing it such that for any element in that open subset, the set $\{(\bar x,\bar y)\in Y_{0}:\;(\bar y,\bar z)\in Y_{0}\;\&\;\bar y.\bar z\in V_{0}'\;\&\;(\bar x,\bar y.\bar z)\in Y_{0}\;\&\;(\bar x.\bar y,\bar z)\in Y_{0}\;\&\;\bar x.(\bar y.\bar z)=(\bar x.\bar y).\bar z\}$ is not large in $Y_{0}$, which means that there exists an open subset of $Y_{0}$ such that one of the following statement fails: 
 $(\bar y,\bar z)\in Y_{0},\;\bar y.\bar z\in V_{0}',\;(\bar x,\bar y.\bar z)\in Y_{0},\;(\bar x.\bar y,\bar z)\in Y_{0}$, or if everything else hold, that $\bar x.(\bar y.\bar z)=(\bar x.\bar y).\bar z$.
 \par We may choose such $\bar z\in V_{0}'$ of the form $\bar z^{\n}$ for some $\bar z\in G$. Since $\{\bar a\in G^*:\;(\bar a,\bar z^{\n})\in Y_{0}\}$ is large in $G^*$, we may choose $\bar a\in G^*$ and $\bar b\in G^*$ of the form $\bar x^{\n}$, $\bar y^{\n}$ with $\bar x,\;\bar y\in G$ such that both $(\bar x^{\n},\bar z^{\n})$ and $(\bar y^{\n},\bar z^{\n})\in Y_{0}$. Moreover we choose
 $\bar y^{\n}$ $\L$-generic over $\bar z^{\n}$ which implies in particular that $\bar y^{\n}.\bar z^{\n}$ is well-defined. Then we choose $\bar x^{\n}$ $\L$-generic over $\bar y^{\n}.\bar z^{\n}$ 
Since $\bar z^{\n}\in V_{0}'$, and since the set of tuples of the form $\bar y^{\n}$ which are $\L$-generic over $\bar z^{\n}$ is dense in $G^*$, $\bar y^{\n}.\bar z^{\n}$ and $(\bar x^{\n}.\bar y^{\n}).\bar z^{\n}$ are well-defined. Since $\bar y^{\n}.\bar z^{\n}\in V_{0}'$, we have $\bar x^{\n}.(\bar y^{\n}.\bar z^{\n})$ is well-defined and that associativity holds for these tuples since they come from elements of $G$. \qed
\medskip
\par Let $V_{1}^{''}:=\{\bar z\in V_{0}':\{ \bar x\in V_{1}:\;(\bar x^{-1},\bar z)\in Y_{0}\;\&\;(\bar x,\bar x^{-1}.\bar z)\in Y_{0}\&\;\bar x.(\bar x^{-1}.\bar z)=\bar z\}$ is large in $V_{0}'\}$; this set is definable by Fact \ref{def-large}.
\cl $V_{1}^{''}$ is large in $V_{0}'$. 
\ecl
\pr Suppose not then there would exist an open subset of $V_{0}'$ over which  $\{ \bar x:\;\bar x.(\bar x^{-1}.\bar z)=\bar z\}$ is not large in $V_{0}'$. So we can find an element of the form $\bar z^{\n}$, $\bar z\in G$ in the complement of that set for which there is an $\bar x\in G$ with $\bar x^{\n}.{(\bar x^{\n})}^{-1}.\bar z^{\n}\neq \bar z^{\n}$, a contradiction. \qed
\medskip 
\par Set $V_{2}:=V_{0}^{''}\cap V_{1}^{''}$.
\cl $V_2$ is open, definable and large in $G^*$.
\ecl
\pr It follows from that both $V_{0}^{''}$ and $V_{1}^{''}$ are large and definable in $G^*$.\qed
\medskip
\par It remains to associate with $G$ a $\L$-definable local group $(G_{local}^*,1,^{-1},\cdot)$ with $^{-1}$ defined on an open subset $U$ of $G_{local}^*$ and $\cdot$ on $O$ an open subset of $G_{local}^*\times G_{local}^*$ satisfying the requirements of Definition \ref{local-def}.
\par Let $O:=\{(\bar x,\bar y)\in (V_{2}\times V_{2})\cap Y_{0}:\;\bar x.\bar y\in V_{2}\;\&\;\forall \bar z\in V_{2}\;((\bar x,\bar y)\in Y_{0}\;\&\;(\bar y,\bar z)\in Y_{0}\;\&\;(\bar x,\bar y.\bar z)\in Y_{0}\;\&\;(\bar x.\bar y,\bar z)\in Y_{0}\rightarrow \bar x.(\bar y.\bar z)=(\bar x.\bar y).\bar z)\}$. It is large in $V_{0}'\times V_{0}'$ and the tuples of the form $(\bar x^{\n},{(\bar x^{\n})}^{-1})$ with $\bar x\in G$ and more generally $(\bar x^{\n},\bar y^{\n})$ with $\bar x, \bar y\in G$ are dense in $O$.
\par The set $U:=\{\bar x\in V_{0}':\;(\bar x,\bar x^{-1})\in O\;\&\;(\bar x^{-1},\bar x)\in O\;\&\;\bar x.\bar x^{-1}=\bar x^{-1}.\bar x=1\,\}$ is large in $V_{0}'$. Suppose not, then we could find $\bar y\in G$ such that $\bar y^{\n}\in V_{0}'\setminus U$ and $(\bar y^{\n},{(\bar y^{-1})}^{\n})\notin O$ (the inverse $^{-1}$ is continuous on $V_{0}'$ and multiplication is continuous on $Y_{0}$), which is a contradiction.
\par Finally we check that $1^{\n}\in U$, that $1^{\n}\times V_{2}\subset O$ and $V_{2}\times 1^{\n}\subset O$. This follows from the fact that $U\subset V_{2}\subset V_{0}'\subset V_{1}$. Since $G^{\n}\cap V_{1}$ was dense in $V_{1}$, the same property holds for $V_{2}$.
\par Then $G^*_{local}:=(V_{2},1,^{-1},\cdot)$ with associated (large) open sets: $U$ and $O$ is a definable local group with the maps $^{-1}$ and $\cdot$ given respectively by the formulas $\phi^*_{-1}$ and $\phi^*_{\times}$.
\qed
\medskip
\rem Recall that one can endow the field $\IR$ of real numbers with a derivation $D$ such that the expansion $(\IR,D)$ is a model of CODF \cite[Theorem 2.5.3]{Br}. Given a definable group $G$ in that expansion $(\IR,D)$, one can associate with $G$ a local definable group $G_{local}^*$ in $\IR$. Then $G_{local}^*$ is necessarily locally Euclidean and so some restriction of $G_{local}^*$ is a local Lie group (moreover some restriction of $G_{local}^*$ is globalisable, namely can be extended to a Lie group).
\erem
\subsection{Subgroups}
\dfn {\rm \cite[Definition 2.10]{G}}\label{sublocal} Let $(H,1,i,m)$ be a local group with associated open sets $U,\;O$. Then a {\it sublocal} group $H_0$ of $H$ (with associated neighbourhood $V$) is a subset $H_{0}$ of $H$ containing $1$ for which there exists an open neighbourhood $V$ of $1$ in $H$ such that:
\begin{enumerate}
\item $H_{0}\subset V$ and $H_{0}$ closed in $V$,
\item for all $x\in H_{0}\cap U$, and $i(x)\in V$, then $i(x)\in H_{0}$,
\item for all $(x,y)\in (H_{0}\times H_{0})\cap O$ and $m(x,y)\in V$, then  $m(x,y)\in H_{0}$.
\end{enumerate}
A sublocal group $H_{0}$ of a definable local group $H$ (with associated neighbourhood $V$) is {\it definable} if $H_{0}$ is definable and $V$ is definable. 
\par For ease of notation we will denote in the following $i$ by $^{-1}$ and $m$ by $\cdot$.
The sublocal group $H_{0}$ of $H$ is {\it normal} if $V$ is symmetric, namely $V\subset U$ and $V=V^{-1}$, and
for all $y\in V$ and $x\in H_{0}$ such that $y\cdot x\cdot y^{-1}$ is defined and belongs to $V$, then $y\cdot x\cdot y^{-1}\in H_{0}.$ 
\edfn
\par Given a (definable) local group $H$ with associated (definable) neighbourhood $U,\;O$ and a sublocal group $H_{0}$ with associated (definable) neighbourhood $V$, one would like to associate with that data a (definable) equivalence relation. Because of the lack of associativity, one has first to restrict oneself to a well-chosen neighbourhood of the identity. We briefly recall below how one proceeds \cite[Section 2]{G}. One first defines the notion for an element $b\in H$ to be represented by a $n$-tuple of elements $(a_{1}, \cdots, a_{n})$ of $H$. Roughly it means that we can write $b$ as a product $a_{1}\cdot \ldots \cdot a_{n}$ and the product can be performed using all the possible choices of parenthesis and gives the same element \cite[Definition 2.4]{G}. Then one says {\it the product $a_{1}\cdot \ldots \cdot a_{n}$ is defined} if there is an element $b\in H$ such that $(a_{1},\cdots,a_{n})$ represents $b$.
Set $\U_{n}^{\times n}:=\U_{n}\times\cdots\times\U_{n}$ with $n$ factors. Then, there is a sequence of decreasing open symmetric neighbourhoods $\U_{n}$ of $1$, $n\in \IN^*$, such that 
for all $(a_{1},\cdots,a_{n})\in \U_{n}^{\times n}$, the product $a_{1}\cdot\ldots\cdot a_{n}$ is defined \cite[Lemma 2.5]{G}. For instance, one may choose $\U_{1}$ as $U\cap U^{-1}$ (so it is definable whenever $U$ is definable),  $\U_{2}$ as an open symmetric neighbourhood of $1$ such that $\U_{2}\times \U_{2}\subset O$ (one considers the open set $(\U_{1}\times \U_{1})\cap O$ and since the topology is definable and is the product topology, we can find an open definable set of the form $\U_{2}\times \U_{2}$ included in that intersection) and in case $\U_{2}$ is not symmetric, we  consider $\U_{2}\cap \U_{2}^{-1}$.
Then let $m_{n}:\U_{2}^{\times n}\rightarrow G:(a_{1}, \cdots, a_{n})\rightarrow a_{1}\cdot \ldots\cdot a_{n}$, one defines inductively in a similar way $\U_{n+1}$ such that $\U_{n+1}^{\times n} \subset m_{n}^{-1}(\U_{2})$. 
\par Let $A\subset \U_{n}$, we will use the notation $A^n:=\{a_{1}\cdot \ldots \cdot a_{n}:\;(a_{1},\cdots,a_{n})\in A^{\times n}\}$.
Now assume that $H$ and $H_{0}$ are definable. Then one can associate a definable equivalence relation with $H_{0}$ (defined in a neighbourhood of the identity) as follows. Choose $W$ a definable symmetric open neighbourhood of $1$ in $H$ such that $W\subset \U_{6}$ and $W^6\subset V$ (for instance choose $W'$ an open definable neighbourhood of $1$ such that $W'^{\times 6}\subset m_{6}^{-1}(V)$ and set $W:=\U_{6}\cap W'\cap W'^{-1}$).
\par Then the binary relation $E_{H_{0}}$ on $W$ defined by $E_{H_{0}}(x,y)$ iff $x^{-1}.y\in H_{0}$ is a definable equivalence relation on $W$ \cite[Lemma 2.13]{G}. Note that Goldbring states the result assuming that $H_{0}$ normal but that hypothesis is only used later when he is dealing with coset spaces (see also \cite[section 1.21]{MZ}).
\par For $x\in W$, define $x.H_{0}:=\{x.h:\;h\in H_{0}\;\&\;(x,h)\in O\}$. Then for $x,\;y\in W$, we have $E_{H_{0}}(x,y)$ iff
$(x.H_{0})\cap W=(y.H_{0})\cap W$.
\par In case $H_{0}$ is in addition normal, this will give rise to a well-defined notion of quotient $H/H_{0}$ which will also be a local group \cite[Lemma 2.14]{G}. 
\medskip
\par Let $G_{0}$ be a $\L_{D}$-definable subgroup of $G$ with $G:=\phi(M)$ and $G_{0}:=\psi(M)$, with $\phi,\;\psi$ two $\L_{D}$-formulas possibly with parameters. The first difficulty is that the corresponding $\L$-formulas $\phi^*,\;\psi^*$ may not have the same number of free variables. For instance, if one takes $G=(M,+,0)$ (with $\phi(x):=(x=x)$) and $G_{0}=(C_{M},+,0)$, where $C_{M}$ is the subgroup of constant elements, (with $\psi(x):=D(x)=0$). Then one can take $\phi^*=\phi$ whereas $\psi^*(x_0,x_{1}):=(x_{1}=0)$ and in this case the associated local groups might be: $G_{local}^*:=(M,+,0)$ and ${G_{0}}_{local}^*:=(\{(u,0):\;u\in M\},+,0)$. In this case in order to view ${G_{0}}_{local}^*$ as a sublocal subgroup of $G^*_{local}$, one has to associate with $G$ the local group $M\times M$.
(A similar example is $G=M\times C_{M}$ and $G_{0}=C_{M}\times C_{M}$; in this case $G^*\subset M^3$ and $G_{0}^*\subset M^4$.)
However, if $G:=\phi(M)$ with $\phi$ with only one free variable (for simplicity) and if whenever $\phi(u)$ holds, it implies that for some $n\geq 0$, the $(n+1)^{th}$ derivative of such element $u$ is in the definable closure of  $u, u^{(1)},\cdots, u^{(n)}$, we may assume that for $G_{0}:=\psi(M)$ a definable subgroup of $G$, then $\psi^*(M)\subset \phi(M)$. If $\phi$ has several free variables, one argues for each free variable separately.
Note that if $G^*$ were definably connected, then $G^*\times M$ is still definably connected.
\lem \label{sublocalgroup} Let $(G,\;G_{0})$ be a  pair of $\L_{D}$-definable groups and assume that $G_{0}$ is a subgroup of $G$. Then one can associate a pair $(H, H_{0})$ of  $\L$-definable local groups to respectively $(G, G_{0})$ in such a way that $H_{0}$ is an $\L$-definable sublocal subgroup of $H$.
In case $G_{0}\trianglelefteq G$, there is a pair of $\L$-definable local groups $(H,H_{0})$ such that $H_{0}$ is normal sublocal subgroup of $H$.
\elem
\pr  Let $\phi,\;\psi$ be $\L_{D}$-formulas defining respectively $G$ and $G_{0}$. By the above discussion we may assume that the corresponding $\L$-formulas $\phi^*$ and $\psi^*$ have the same number of free variables. We proceed as in Theorem \ref{local} to associate with $G$ an $\L$-definable local group $G^*_{local}$ with $V_2$ the domain of $G^*_{local}$ and data $(U,O,i,m)$. Let $\tilde V_{2}:=\{\bar x\in V_{2}:\ (\psi^*_{mod}(\bar x)\;\&\;\bar x\in U)\rightarrow \psi^*_{mod}(i(x))\}$. This is a large subset of $V_{2}$. If not we could find $\bar x\in G_{0}$ and so $\psi^*_{mod}(\bar x^{\n})$ and $\psi^*_{mod}(i(\bar x^{\n}))$, a contradiction since $\bar x^{-1}\in G_{0}$.
Set $\tilde O:=\{(\bar x,\bar y)\in O:\,(\psi^*_{mod}(\bar x)\;\&\;\psi^*_{mod}(\bar y)\rightarrow \psi^*_{mod}(m(\bar x,\bar y))\}$. Again this is a large subset of $O$. Let $H_{0}$ be the relative closure of $\psi^*_{mod}(M)\cap \tilde V_{2}$ in $\tilde V_{2}$. Define $H:=
\tilde G^*_{local}$ a new local group associated with $G$ with domain $V_{2}$ and data $(U,\tilde O,i,m)$. Then $H_{0}$ is a sublocal subgroup of $H$ with associated neighbourhood  $\tilde V_{2}$. One easily checks the conditions (1) up to (3) of Definition \ref{sublocal} (it follows from our construction and the fact that the closure of a subgroup in  a topological group is still a subgroup.)
\par Now assume that $G_{0}\trianglelefteq G$ and let $H$, respectively $H_{0}$ the local groups associated with $G$, respectively $G_{0}$ and let $V$ be associated with $H_{0}$. We define $\tilde V\subset V$ such that the restriction of $H_{0}$ to $\tilde V$ is normal in $H$.
\par Consider $\tilde V:=\{\bar y\in V:\{\bar z\in G_{0}^*:\bar y.\bar z.\bar y^{-1}$ is defined and belongs to $G_{0}^*\}$ is large in ${G_{0}}^*\}$. Let us show that $\tilde V$ is large in $V$. Suppose not then there is a tuple $\bar y^{\n}\in V\setminus \tilde V$ (in particular $\bar y\in G$). So, the set $\{\bar z\in G_{0}^*:\bar y.\bar z.\bar y^{-1}$ is defined and belongs to $G_{0}^*\}$ is not large in ${G_{0}}^*$. So there is a tuple $\bar z^{\n}$ such that $\bar z\in G_{0}$ and $\bar y^{\n}.\bar z^{\n}.\bar y^{{\n}^{-1}}\notin G_{0}^*$, a contradiction with the fact that $G_{0}\trianglelefteq G$.
\qed
\medskip
\par So given two $\L_{D}$-definable groups $G_0,\;G$ with $G_{0}$ a subgroup of $G$, and $H$, respectively $H_{0}$ the local groups associated with $G$, respectively $G_{0}$ and $V$ associated with $H_{0}$ in such a way that $H_{0}$ is a sublocal subgroup of $H$. On a definable neighbourhood $W$ of $1$ (defined as above), we get a definable equivalence relation $E_{H_{0}}$. In case $G_{0}$ is a normal subgroup of $G$ and the quotient $G/G_{0}$ is $\L_{D}$-definable, is the local group associated with the quotient $H/H_{0}$ equivalent (in the sense of \cite[Lemma 2.15]{G}) to the local group associated to $G/G_{0}$? 
\medskip
\thm \label{iddc} Let $\M\models CODF$ and $(G_{j})_{j\in \omega}$ be a descending chain of $\L_{D}$-definable groups in $\M$. Then for finitely many $j$, we have that the index $[G_{j}:G_{j+1}]$ is infinite and the number of such $j's$ is bounded by $\L$-dim$({G_{0}}^*_{local})$.
\ethm
\pr Let $\M\models CODF$ and $(G_{j})_{j\in \omega}$ be a descending chain of $\L_{D}$-groups in $\M$ and let $\phi_{j}(x_{1},\cdots,x_{n})$ be an $\L_{D}$-formula with $\phi_{j}(M)=G_{j}$. For $1\leq i\leq n$, either there exists $j(i)$ such that $\phi_{j(i)}$ implies a condition on the successive derivatives of $x_{i}$ and so there exists $n(i)\geq 1$ such that if $\phi_{j(i)}(u_{1},\cdots, u_{n})$ holds, then $u_{i}^{(n(i))}\in dcl(u_{i},\cdots,u_{i}^{(n(i)-1)})$, or there is no such index $j(i)$ and set $n(i)=0$. Let $k=max_{1\leq i\leq n} n(i)+1$. Then we consider $G_{0}^*$ as a subset of $M^n.k$ and we may assume that for each $j\leq \ell\in \omega$ that $G_{j}^*\subset G_{\ell}^*$.
\par By Lemma \ref{sublocalgroup}, given $G_{j+1}\subset G_{j}$, we have associated local groups $H_{j+1},\;H_{j}$ in such a way $H_{j+1}$ is a sublocal subgroup of $H_{j}$. We begin with $G_{0}$ with associated local group $H_{0}$ and neighbourhood $(U,O,i,m)$, then given $G_{j+1}\subset G_{j}$, $j\geq 0$, we get a sublocal subgroup $H_{j+1}$ of $H_{j}$ with associated neighbourhood $V_{j+1}$ which is a large subset of $V_{j}$. Recall that going from $G_{0}$ to $G_{1}$, we modified $O$ to a large subset $O_{1}$ of $O$ and so we repeat the procedure going from $G_{j}$ to $G_{j+1}$. Let us call the modified open large subset of $O$, $O_{j+1}$ and the corresponding neighbourhood of $1$ where $H_{j+1}$ induces an equivalence relation, $W_{j}$. Note that if $[G_{0}:G_{j}]$ is finite, then $\L$-dim$(H_{0})$=$\L$-dim$(H_{j})$.
\par Let $j\geq 0$ be the smallest index such that $[G_{j}:G_{j+1}]$ is infinite. 
Let us show that  $\L$-dim$(H_{j})<\L$-dim$(H_{j+1})$. Suppose otherwise, then the associated equivalence relation $E_{H_{j+1}}$ induces distinct equivalence classes on $H_{j}$ of dimension equal to $\L$-dim$(H_{j+1})$. So they have non-empty interior in $H_{j}$ and so we have finitely many distinct such equivalence classes \cite[Proposition 2.1]{P86}. Therefore the index $[G_{j}:G_{j+1}]$ is finite, a contradiction.
\par Since $\L$-dim$(H_{0})$ is finite, this can happen only finitely many times.
\par Finally note that if $\L$-dim$(H_{0})=0$, then $H_{0} $ is finite and so $G_{0}$ is finite.
\qed
\medskip
\par Recall that an $M_{c}$-group is a group satisfying the descending chain condition on centralisers
\cor Let $\M\models CODF$ and let $G$ be a definable subgroup in $\M$. Then $G$ is an $M_{c}$-group. In particular, the Fitting subgroup $F(G)$ of $G$ is nilpotent and definable.
\ecor
\pr Let $(G_{j})_{j\in \omega}$ be a chain of centralisers of the form $G_{j}:=C(A_{j})$ where $A_{j}\subset G$, $A_{j}$ finite and $A_{j+1}\supset A_{j}$. 
In this case if $x\in C(A_{j})$, then $x^{-1}\in C(A_{j})$ and if $x,\;y\in C(A_{j})$, then $x.y\in C(A_{j})$.
Indeed $x.a=a.x$ implies that $(x.a)^{-1}=(a.x)^{-1}$, so $a^{-1}.x^{1}=x^{-1}.a^{-1}$ and multiplying by $a$ on the left and on the right, we get
$a.(a^{-1}.x^{1}).a=a.(x^{-1}.a^{-1}).a$ and so using associativity, and $a.a^{-1}=a^{-1}.a=1$, we get $x^{-1}\in C(a)$.
To show that $x,\;y\in C(a)$ implies that $x.y\in C(a)$, only uses associativity. As before let $H_j$ be the local group associated with $G_j$ with corresponding neigbourhoods $V_j$ and $W_j$. So in this case when going from $H_{j}$ to $H_{j+1}$, we can keep the same associated neighbourhood $V_{j}$ and the same neighbourhood $W_{j}$ on which $H_{j+1}$ induces an equivalence relation, namely we do not need to restrict them at each step. So we apply \cite[Proposition 2.1]{P86} to show that there are only finitely many subgroups $G_{j}$ of finite index. We conclude using Theorem \ref{iddc}.
\par Then it suffices to apply a former result of F. Wagner \cite[Theorem 1.2.11]{W} to get that the Fitting subgroup $F(G)$ is nilpotent. To see that $F(G)$ is definable, one either uses an observation of Ould Houcine, or \cite[Fact 3.5]{BJO}. \qed

\subsection{Type-definable groups}
Throughout this subsection, we will assume that $\kappa$ is a strongly inaccessible cardinal and that $\M$ is a $\kappa$-saturated model of $T_{c,D}^*$. We will say that a set is {\it bounded} if it is of cardinality strictly less than $\kappa$. A type-definable subset of $M^n$ is a set defined by a bounded set of formulas.
\par Let $G$ a definable group in $\M$ and w.l.o.g. we may assume that $G$ is $\emptyset$-definable. Let $A\subset M$ and denote by $G_{A}^{00}$ the minimal $A$-type definable subgroup of bounded index. Note that this subgroup is normal \cite[Remark 2.9]{P04}. One says that $G^{00}$ exists if $G^{00}_{A}=G^{00}_{\emptyset}$ for all $A$ \cite[Definition 1.7]{KS}. S. Shelah proved that for any definable group $G$ in a NIP theory, $G^{00}$ exists \cite[Theorem 1.8]{KS}, \cite{Sh}. One can endow the quotient with the so-called logic topology in such a way $G/G^{00}$ become a compact topological group \cite[Lemma 2.5]{P04}.
\par One also defines $G^0$ as the intersection of all (relatively) definable subgroups of finite index \cite[8.1]{simon}, it is again a normal $\emptyset$-type-definable subgroup of bounded index. 
\medskip
\par In Lemma below, we will use the notation introduced in Proposition \ref{*densek} and Theorem \ref{local}.
\lem \label{open}
Assume that $G$ is a definable subgroup of $\M$
and let $H$ be a definable subgroup of $G$ and assume that $H^*\subset G^*$. Let $G_{local}^{*}$ be the local group associated with $G$ with data $(U,O,{\;}^{-1},\cdot)$. Assume that $H^*\cap G_{local}^{*}$ has non empty relative interior in $G_{local}^*$, then
$H^{*}\cap U$ is a relatively open subset of $G_{local}^{*}$.
\elem
\pr Let $G_{local}^*$ a local group associated with $G$ with domain $V_2$ and data $(U, O, {\,}^{-1},\cdot)$.
Since $H^*\cap G_{local}^{*}$ has non empty relative interior in $G_{local}^*$, there is an open definable subset $\tilde U$ of some cartesian product of $M$ such that $\tilde U\cap G_{local}^*\subset H^*\cap G_{local}^{*}$. W.l.o.g. we may assume that $\tilde U\subset U$ and that $1^{\n}\in \tilde U$. (Otherwise, we choose an element of the form $h^{\n}\in H^*\cap \tilde U$ and by Theorem \ref{local}, multiplication by $(h^{\n})^{-1}$ is defined and continuous on $U$. So we consider $(h^{\n})^{-1}.\tilde U$. Since $h^{\n}\in U$, $(h^{\n},(h^{\n})^{-1})\in O$ and so $1^{\n}\in (h^{\n})^{-1}.\tilde U$).
\par Let $\bar h\in H^*\cap U$ (therefore $\bar h^{-1}$ is well-defined). The inverse function is continuous on $U$ and so given $U_{1}\subset U$ a neighbourhood of $1^{\n}$ such that $U_{1}\cap G_{local}^*\subset H^*\cap G_{local}^{*}$, we have that there exists $U_{\bar h}\subset U$ a neighbourhood of $\bar h$ such that $\bar h^{-1}.U_{\bar h}\subset U_{1}$. Moreover we may assume that $U_{\bar h}\times U_{1}\subset O$. Let us show that $H^*\cap U=\bigcup_{\bar h\in H^*\cap U} U_{\bar h}$. 
\par We have that $U_{\bar h}\subset \bar h.U_{1}$. Indeed, $(\bar h^{-1}.\bar w)\in U_{1}$ with $\bar w\in U_{\bar h}$ and $\bar h.(\bar h^{-1}.\bar w)=\bar w$, since $(\bar h, \bar h^{-1}.\bar w)\in O$ and $(1^{\n},\bar w)\in O$. 
\qed
\dfn Let $X$ be a type-definable subset of $M^k$. Then we denote by $X^*$ the subset $\bigcap_{\ell\in L} X_{\ell}^*$, where $\vert L\vert<\kappa$ and $X=\bigcap_{\ell\in L} X_{\ell}$, $X_{\ell}$ definable, where $X_{\ell}^*$ has been defined in Proposition \ref{*densek}. We will always assume that this family $X_{\ell},\;\ell\in L$ is closed under finite intersection. Note that there is always $n\in \IN$ such that $X^*\in dcl(\tilde X^*)$, where $\tilde X^*\subset M^n$. 
Indeed if $(x_{1},\cdots,x_{k})\in X$, then for every component $x_{j}$ of $\bar x$, $1\leq j\leq k$, there is $\ell$ such that the condition that $\bar x\in X_{\ell}$ imposes a condition on the $i^{th}$-derivative of $x_{j}$ for some $i\geq 1$, or for no $\ell$ there is a  condition on the successive derivatives of $x_{j}$. Therefore we may assume that in the presentation $(\star)$ of $X$, we have the property that each $X_{\ell}^*\subset M^{n_{1}}\times\cdots\times M^{n_{k}}$, $\L$-dim$(X_{\ell}^*)=m=\sum_{i=1}^k m_{i}$, with $m_{i}\leq n_{i}$, $1\leq i\leq k$ and if $X_{\ell}\subset X_{\ell'}$, then $X_{\ell}^*\subset X_{\ell'}^*.$
\edfn
With the same conventions as in the above Lemma, we get:
\lem Let $X$ be a type-definable subset of $M^k$, $X=\bigcap_{\ell\in L} X_{\ell}$, where $X_\ell$ is an $\L_D$-definable subset of $M^k$. Then $X^*$ is relatively open in each $X_{\ell}^*$, $\ell\in L$.
\elem
\pr Let $n$ be such that each $X_{\ell}^*\subset M^n$ and let $m:=\L$-dim$(X_{\ell}^*)$. In each $X_{\ell}^*$, choose a generic point of the form $u_{\ell}^{\n}$. 
Let us show that each  $X_{\ell}^{*}$ has non-empty interior in $X^{*}$. 
\par Since $u_{\ell}^{\n}$ is a generic point of $X_{\ell}^*$, there is a subset $U_{\ell}$ containing $u_{\ell}^{\n}$ and a projection $\pi$ into $M^m$ such that $\pi(U_\ell)$ is an open subset of $\pi(X_{\ell}^*)$ with $U_\ell\subset dcl(\pi(U_\ell))$.
Recall that a basis of neigbourhoods of zero (in $M$) is given by $V_{\bar y}:=\{\varrho(M,\bar y):\;\bar y\in M\}$.  
So we may assume that for some $\bar b_{\ell}\subset M$, $U_\ell$=$u_{\ell}^{\n}+\widehat{V_{\bar b_{\ell}}^{m}}$, where $\widehat{V_{\bar b_{\ell}}^{m}}\subset dcl(V_{\bar b_{\ell}}^{m})$ and is $\L$-definable.  
We write the following type $p(x,\bar y):=\{x^{\n}\in X_{\ell}^*\,\&\, (x^{\n}+\widehat{V_{\bar y}^m})\subset X_{\ell}^*, \;\ell\in L\}$. 
This is finitely satisfiable (since we assumed that the family $X_{\ell}$ was closed by finite intersection, it suffices to check it for each $\ell$, so take $\bar x=u_{\ell}^{\n}$ and $\bar y=\bar b_{\ell}$ such that $U_\ell$=$u_{\ell}^{\n}+\widehat{V_{\bar b_{\ell}}^{m}}$). Since this is a small set of formulas, it is satisfiable in $\M$. Let $(a,\bar b)$ be a realization of that type. Then $a^{\n}\in X$ and $a^{\n}+\widehat{V_{\bar b}^m}\subset X_{\ell}^*$, for every $\ell\in L$. \qed

\prop \label{type-open} Let $G$ be a definable group of $\M$ and $G^*_{local}$ be a $\L$-definable local group associated with $G$ with data $(U,O,{\;}^{-1},\cdot)$.
Let $H$ a type-definable subgroup of bounded index in $G$, then $H^*\cap U$ is a relatively open subset of $G_{local}^*$. In particular, when $T_{c}$ is in addition NIP, both ${(G^0)}^*$ and ${(G^{00})}^*$ are of the same topological dimension and this dimension is equal to $\L$-dim$(G^*)$.
\eprop
\pr W.l.o.g. $H$ is a bounded intersection of definable subsets $X_{\ell}$, $\ell\in L$, closed under finite intersection. By compactness, for each $\ell\in L$, finitely many translates of $X_{\ell}$ cover $G$ (otherwise neither a bounded number). So each $X_{\ell}$ has the same t-dimension as $G$.
\par For each $\ell\in L$, let $g_{\ell i}\in G$, $i\in I_{\ell}$ be such that $\bigcup_{i\in I_{\ell}} g_{\ell i}X_{\ell}=G$. 
Assume that $G$ is defined by an $\L_{D}$-formula $\phi(\bar x)$ and $X_{\ell}$ by an $\L_{D}$-formula $\psi_{\ell}(\bar x)$. Consider the following subset $Z\subset G_{local}^*$ defined by $Z:=\{\bar x\in G^*_{local}:\;\phi^*(\bar x)\leftrightarrow \bigvee_{i\in I_{\ell}} \exists \bar y\;(\psi_{\ell}^*(\bar y)\;\&\;(g_{\ell i}^{\n},\bar y)\in O\;\&\;\bar x=g_{\ell i}^{\n}.\bar y)\}$. Since $Z$ contains $G^{\n}\cap G^*_{local}$, $Z$ is large in $G^*_{local}$ and so of the same $\L$-dimension. Multiplication by $g_{\ell i}^{\n}$ is a bijection when defined on $G^*_{local}$ and so each $X_{\ell}^*:=\psi_{\ell}^*(M)$ has the same $\L$-dimension as $G_{local}^*$. This implies that each $X_{\ell}^*$ has non-empty interior in $G^*_{local}$. 
\par By the above Lemma, $H^*$ has non-empty interior in each $X_{\ell}^*$, which implies by the above that $H^*$ has non-empty interior in $G^*_{local}$. Since $H$ is a subgroup, $H^*\cap U$ is a relatively open subset of $G_{local }^*$ (we can use the last part of the argument of  Lemma \ref{open} which did not use that $H$ was definable).
\qed
\medskip
\par In Theorem \ref{local}, we show how to associate with a $\L_{D}$-definable group $G$ in a model of $T_{c,D}^*$, a $\L$-definable local group $G^*_{local}$ with associated data $(U,O,i,m)$ such that $G^{\n}\cap G^*_{local}$ is dense in $G^*_{local}$. When $T_{c}$ is NIP, we know that $G^{00}$ exists and that $G/G^{00}$ with the logic topology is a compact topological group. We would like to put on $G^*_{local}$ (or more precisely on a large open subset) a bounded type-definable equivalence relation $\sim^{00}$ associated with $(G^{00})^*$. 
\prop \label{equiv} Assume (in addition) that $T_{c}$ is NIP and let $G$ be a $\L_{D}$-definable group of $\M$. Then there is a bounded type definable relation $\sim^{00}$ on $G^*_{local}$, which is an equivalence relation on a large subset of $G^*_{local}$. Moreover, for any $a,\;b\in G$ with $a.b^{-1}\in G^{00}$, we have that $a^{\n}\sim^{00} b^{\n}$. 
\eprop
\pr As in the preceding proposition, we write $G^{00}$ as a bounded intersection of definable subsets $X_{\ell}$, $\ell\in L$, closed under finite intersection and with the property that for each $\ell\in L$, there is a finite set of indices $I_{\ell i}$ and $g_{\ell i}\in G$, $i\in I_{\ell}$ such that $\bigcup_{i\in I_{\ell}} g_{\ell i}X_{\ell}=G$.
Moreover, by Proposition \ref{type-open} and its proof, each $X_{\ell}^*$ has non-empty interior in $G^*_{local}$ and $\bigcup_{i\in I_{\ell}} g_{\ell i}^{\n}.X_{\ell}^*$ is a large subset of $G^*_{local}$.
\par So we define on $(U\times U)$ a (symmetric) relation $\sim$ as follows:  $$\bar x\sim \bar y\;\;{\rm iff }\;\;(\bar x,\bar y^{-1})\in O\;\&\;(\bar y,\bar x^{-1})\in O\;\&\;\bigwedge_{\ell\in L}\;\bar x.\bar y^{-1}\in X_{\ell}^*\;\&\;\bar y.\bar x^{-1}\in X_{\ell}^*.$$ 
Then we define $$\bar x\sim^{00} \bar y\;\;{\rm  iff}\;\; \exists x_{0}\in G\exists y_{0}\in G\;(\bar x\sim x_{0}^{\n}\;\&\;x_{0}^{\n}\sim y_{0}^{\n}\;\&\;y_{0}^{\n}\sim \bar y).$$
We consider the following subset of $G^*_{local}$: $\{\bar x\in G^*_{local}\cap U\cap \U_6:\forall x_{0}\in G\;
\forall x_{1}\in G\,((\bar x.{(\bar x_{0}^{\n})}^{-1}\in {G^{00}}^*\;\&\;\bar x_{0}^{\n}.\bar x^{-1}\in {G^{00}}^*\;\&\;
\bar x_{1}^{\n}.\bar x^{-1}\in {G^{00}}^*\;\&\;\bar x.{({\bar x_{1}}^{\n})}^{-1}\in {G^{00}}^*)\rightarrow (\bar x_{1}.\bar x_{0}^{-1}\in G^{00}\;\&\;\bar x_{0}.\bar x_{1}^{-1}\in G^{00}))\}.$



This is a large subset of $G^{*}_{local}\cap U\cap \U_6$ that we will denote it by $\widehat{G^*_{local}}$. Otherwise we could find an element of the form $\bar u^{\n}$ with $\bar u\in G$ in the complement and so there would exist $\bar x_{0},\bar x_{1}\in G$ such that $(\bar u^{\n}.{(\bar x_{0}^{\n})}^{-1}\in {G^{00}}^*\;\&\;
\bar x_{1}^{\n}.{(\bar u^{\n})}^{-1}\in {G^{00}}^*)$, but $\bar x_{1}.\bar x_{0}^{-1}\notin G^{00}$, a contradiction since $G^{00}$ is a subgroup.
On that large subset, let us check that $\sim^{00}$ is an equivalence relation. 
This follows from the fact that if $\bar x\sim \bar x_{0}^{\n}$ and $\bar x\sim \bar x_{1}^{\n}$, then $\bar x_{0}^{\n}\sim \bar x_{1}^{\n}$.
\medskip
\par Moreover, if $a, b\in G$, then ($a^{\n}\sim^{00} b^{\n}$ iff for all $\ell\in L$, $a^{\n}.(b^{\n})^{-1}\in X_{\ell}^*$), equivalently iff $a.b^{-1}\in G^{00}$ (note that $a^{\n}.(b^{\n})^{-1}=(a.b^{-1})^{\n}$).

\qed
\medskip
\par One can endow the quotient $\widehat{G^*_{local}}/\sim^{00}$ with the following {\it logic topology} \cite{LP}. Let $\mu:G^*_{local}\rightarrow G^*_{local}/\sim^{00}$ be the canonical surjection and define the closed subsets $Z$ of $\widehat{G^*_{local}}/\sim^{00}$ as those such that there exists $Y$ an $\emptyset$-type-definable subset of $\M$ such that $\mu(Y)\cap \widehat{G^*_{local}}=Z$ (see \cite[Definition 2.3]{P04}). 
\cor The quotient $\widehat{G^*_{local}}/\sim{00}$ endowed with the logic topology is a Hausdorff compact topological space.
\ecor
\pr This follows from Propositions \ref{type-open}, \ref{equiv} and \cite[Lemma 2.5]{P04}. We use the fact that a small intersection of formulas defines an open subset of $G^*_{local}$. \qed
\medskip
\par Let $G^{00}=\bigcap_{\ell\in L} X_{\ell}$; we know that $G^{00}$ is a normal subgroup of $G$ and so letting $g_{i}\in G$, $i\in I$, be coset representatives of $G^{00}$ in $G$, with $\vert I\vert<\kappa$, we have that $G^{00}=\bigcap_{\ell\in L,\;i\in I} X_{\ell}^{g_{i}}$ with $X_{\ell}^{g_{i}}=g_{i}.X_{\ell}.g_{i}^{-1}$.
\medskip
\par A natural question is to see whether with the logic topology $\widehat{G^*_{local}}/\sim^{00}$ is a local group and under which conditions it becomes a local Lie group.
\par I. Goldbring has shown that every locally compact NSS local group has a restriction which is a local Lie group and that every locally Euclidean local group is NSS \cite{G}.

\medskip
\par Recall that a group is {\it definably connected} if it has no definable subgroup of finite index. As recalled above if $T_c$ is NIP, then given an $\L_D$-definable group in a model $\M$ of $T_{c,D}^*$, $G^0$ is a type-definable, definably connected subgroup of $G$.
\ft \label{connect} {\rm \cite[Fact 1.5 (iii)]{BOPP}} Assume that $G$ is definably connected and $H$ is a normal type-definable subgroup of $G$ of bounded index. Then, $G/H$ is connected.
\eft
\pr For convenience of the reader we give a proof. As in the above we endow the quotient $G/H$ with the logic topology and so it becomes a compact topological group. As such one can define the connected component of the identity, that we will denote by $(G/H)_{0} $. Since $G/H$ is compact, $(G/H)_{0}$ is the intersection of all normal open subgroups of $G/H$ \cite[chapter 2, Theorem 4]{Hig}. If $G/H$ were not connected, one of these subgroups, say $N$, would be proper and since $G/H$ is compact, $N$ is a clopen proper subgroup of finite index. By definition of the topology we put on $G/H$, the inverse image of $N$ is a type-definable subgroup of $G$ of finite index. So this inverse image is a proper definable subgroup of finite index in $G$, a contradiction.
\qed
\cor Assume that $G$ is definably connected and that $H$ is a normal type-definable of bounded index, then either $G/H$ is abelian, or it contains a free group of rank $2^{\aleph_{0}}.$ 
\ecor
\pr Since $G$ is definably connected, $G/H$ is connected by Fact \ref{connect}. Since $G/H$ is a compact group, we may apply a theorem of Balcerzyk and Mycielski: either $G/H$ is abelian, or it contains a free group of rank $2^{\aleph_{0}}$ \cite[Theorem 2]{BM}.
\qed

\bigskip

\par {\bf Acknowledgments:} I would like to thank  the {\it International Center Mathematical Science (ICMS, Edinburgh)}, where I had the opportunity to give a talk on that subject in july 2014.

\end{document}